\documentclass[12pt]{amsart} 
\usepackage{epsf}

\def\ie{{i.e.}}

\def\viceversa{{\it vice versa}}
\def\cf{{\it cf.}}
\def\fg{fundamental group}
\def\fgs{fundamental groups}
\def\wid{{wiring diagram}} 

\newcommand\set[1]{{\{{#1}\}}}
\def\suchthat{{\,:\,\,}}
\newcommand\power[1]{{|{#1}|}}
\def\Wmodequiv{{W_S/\!\equiv}}

\newcommand\sg[1]{{\left<{#1}\right>}}

\newtheorem{thm}{Theorem}[section]
\newtheorem{prop}[thm]{Proposition}
\newtheorem{cor}[thm]{Corollary}
\newtheorem{lem}[thm]{Lemma}
\newtheorem{clm}[thm]{Claim}
\newtheorem{conj}[thm]{Conjecture}
\newtheorem{exa}[thm]{Example}
\newtheorem{defn}[thm]{Definition}

\newtheorem{rem}[thm]{Remark}

\newcommand{\ben}{\begin{enumerate}}
\newcommand{\een}{\end{enumerate}}
\newcommand{\ble}{\begin{lem}}
\newcommand{\ele}{\end{lem}}
\newcommand{\bre}{\begin{rem}}
\newcommand{\ere}{\end{rem}}
\newcommand{\bth}{\begin{thm}}
\newcommand{\eth}{\end{thm}}
\newcommand{\bpr}{\begin{prop}}
\newcommand{\epr}{\end{prop}}
\newcommand{\bco}{\begin{cor}}
\newcommand{\eco}{\end{cor}}
\newcommand{\bcon}{\begin{conj}}
\newcommand{\econ}{\end{conj}}
\newcommand{\bde}{\begin{defn}}
\newcommand{\ede}{\end{defn}}
\newcommand{\bex}{\begin{exa}}
\newcommand{\eex}{\end{exa}}
\newcommand{\bcl}{\begin{clm}}
\newcommand{\ecl}{\end{clm}}
\newcommand{\barr}{\begin{array}}
\newcommand{\earr}{\end{array}}
\newcommand{\btab}{\begin{tabular}}
\newcommand{\etab}{\end{tabular}}
\newcommand{\beq}{\begin{equation}}
\newcommand{\eeq}{\end{equation}}
\newcommand{\bea}{\begin{eqnarray*}}
\newcommand{\eea}{\end{eqnarray*}}
\newcommand{\bce}{\begin{center}}
\newcommand{\ece}{\end{center}}
\newcommand{\bpi}{\begin{picture}}
\newcommand{\epi}{\end{picture}}
\newcommand{\bfi}{\begin{figure} \begin{center}}
\newcommand{\efi}{\end{center} \end{figure}}

\newcommand{\bsl}{\begin{slide}{}}
\newcommand{\esl}{\end{slide}}

\newcommand{\bib}{thebibliography}

\newcommand{\sbs}{\subset}

\newcommand\Lpair[2]{{\left[{#1},{#2}\right]}}

\newcommand{\ra}{\rightarrow}

\newcommand{\be}{\beta}

\newcommand{\si}{\sigma}

\newcommand\G\Gamma
\newcommand{\Ga}{\Gamma}

\newcommand{\ba}{{\bf a}}
\newcommand{\bb}{{\bf b}}

\newcommand{\bw}{{\bf w}}

\newcommand{\C}{{\mathbb C}}

\newcommand{\R}{{\mathbb R}}

\newcommand{\PP}{{\mathbb P}}
\newcommand{\cal}{\mathcal}

\newcommand{\cB}{{\cal B}}

\newcommand{\cL}{{\cal L}}

\newcommand{\Int}{\mathop {\rm Int}}

\newcommand{\LPs}{Lefschetz pairs}
\newcommand{\arel}{{\stackrel{\Delta}{=}}}

\newcommand\defin[1]{{\it{#1}}}
\newcommand\figs[1]{#1}
\newcommand\FIGUREy[4][]{{\begin{figure}[!h]\epsfysize=#3 {\epsfbox{\figs{#2}}}\caption{#1}\label{#4}\end{figure}}}
\newcommand\FIGUREx[4][]{{\begin{figure}[!h]\epsfxsize=#3 {\epsfbox{\figs{#2}}}\caption{#1}\label{#4}\end{figure}}}

\newcommand\tru[3][c+]{{{#1}\triangle_{#2}^{(#3)}}}
\newcommand\trd[3][c+]{{{#1}\nabla_{#2}^{(#3)}}}
\long\def\forget#1\forgotten{}

\newcommand\eq[1]{{(\ref{#1})}}
\def\Cdots{{\cdot\dots\cdot}}
\renewcommand\th[1]{{$#1$th}}

\begin{document}

\title[Classes of wiring diagrams]%
{Classes of wiring diagrams and their invariants}

\author[David Garber, Mina Teicher]{David Garber$^1$, Mina Teicher$^1$}
\email{\{garber,teicher\}@macs.biu.ac.il}
\address{
Dept. of Math. and CS, Bar-Ilan University, Ramat-Gan 52900,
Israel
}
\author[Uzi Vishne]{Uzi Vishne$^2$}
\email{vishne@math.huji.ac.il}
\address{
Einstein Institute of Mathematics, Hebrew University, Jerusalem
91904, Israel
}

\date{\today}

\stepcounter{footnote}
\footnotetext{Partially supported by The Israel Science Foundation
(Center of Excellence Program), by the Emmy Noether Institute for
Mathematics and by the Minerva Foundation (Germany).
The research was done during the Ph.D. studies of David Garber,
under the supervision of Prof. Mina Teicher.}
\stepcounter{footnote}
\footnotetext{Partially supported by the Edmund Landau Center for
Research in Mathematical Analysis and Related Subjects.}

\begin{abstract}
Wiring diagrams usually serve as a tool in the study of arrangements of lines and
pseudolines. Here
we go in the opposite direction, using known properties of line
arrangements to motivate certain equivalence relations and actions on
sets of wiring diagrams, which preserve the incidence lattice and the fundamental groups
of the affine and projective complements of the diagrams.

These actions are used in \cite{GTV} to classify real arrangements of up to $8$ lines and
show that in this case, the incidence lattice
determines both the affine and the projective fundamental groups.
\end{abstract}

\maketitle


\section{Introduction}

A \defin{line arrangement} in $\C^2$ is a union of finitely many
copies of $\C^1$. An arrangement is called \defin{real } if the
defining equations of its lines can be written with real
coefficients, and \defin{complex} otherwise.

The fundamental groups of the complement in $\C\PP^2$ and of its
affine part in $\C ^2$ are called the \defin{projective} and
\defin{affine \fgs} of the arrangement, respectively.

A more combinatorial invariant is the incidence lattice of the arrangement,
which consists of the lines, their intersection points,
the empty set and $\C\PP ^2$, ordered by inclusion.

To a real arrangement, one can associate a combinatorial object, called a wiring diagram.
This is not an invariant, as it depends
on the choice of a guiding generic line. Still, this object turns out
to be a useful tool in the study of arrangements of
lines, \cf\ \cite{MoTe1}, \cite{GTV} or \cite{CS}.

\medskip

In this paper, we show how the natural actions on the plane induce
actions on sets of wiring diagrams. Likewise, identities in the
braid group induce some equivalence relations on these sets. The
actions and relations will be shown to preserve  the incidence
lattice and the affine and projective fundamental groups. We also
study their relations and show that they induce dihedral group
actions on the set of diagrams.

In general, there are millions of wiring diagrams of up to $8$
lines, each equipped with its own fundamental groups. In
\cite{GTV} we use the actions and relations intorduced here to show
that there are at most $103$ distinct fundamental groups, and
they are determined by the incidence lattice of the diagram.

\medskip

The paper is organized as follows. In section \ref{defs} we
present some combinatorial objects which are associated to real
wiring diagrams: the list of \LPs, the incidence lattice, and the
signature. Section \ref{MT} describes the algorithm commonly used
for computing finite presentations of the affine and projective
fundamental groups.

In section \ref{actions} we define and study  equivalence
relations on the set of wiring diagrams, and discuss certain
actions on the equivalence classes. We prove that these relations
and actions preserve the incidence lattice and the affine and
projective fundamental groups. Finally, in section
\ref{connections} we discuss some connections between the various
actions, and a numerical example. The focus in this paper is on
the actions and relations, and their properties. The application
to classification of arrangements with up to $8$ real lines is
given in \cite{GTV}.

\section{Combinatorial preliminaries} \label{defs}

We briefly recall some of the combinatorial objects and constructions
related to arrangements of lines.

\subsection{Wiring diagrams and Lefschetz pairs}

Given a real arrangement of lines in $\C^2$, the intersection of
its affine part with the natural copy of $\R^2$ in $\C^2$ is an
arrangement of lines in the real plane.

To an arrangement of $\ell$ lines in $\R^2$ one can associate a
\defin{wiring diagram} \cite{Go}, which holds the combinatorial
data of the arrangement and the position of the intersection
points. A wiring diagram is a collection of $\ell$ wires (where a
\defin{wire} in $\R ^2$ is a union of segments and rays,
homeomorphic to $\R$). The induced wiring diagram is constructed
by choosing a new line (called the \defin{guiding line}), which
avoids all the intersection points of the arrangement, such that
the projections of intersection points do not overlap. Then,
$\ell$ wires are generated as follows. Start at the '$\infty$'
end of the line with $\ell$ parallel rays, and for every
projection of an intersection point, make the corresponding
switch in the rays, as in Figure \ref{latowd}.

\FIGUREy{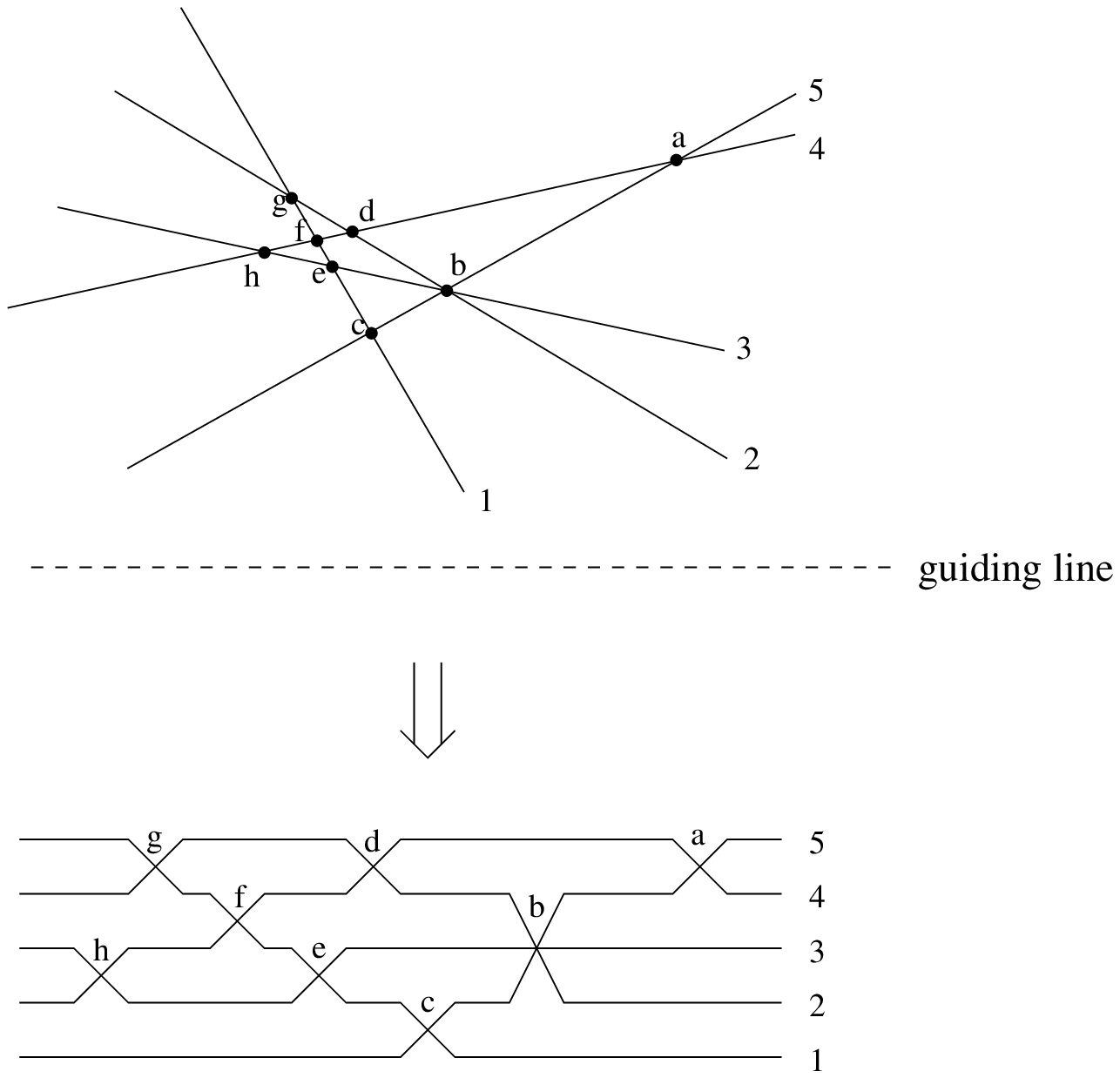}{6cm}{latowd}


To a wiring diagram, one can associate a list of \defin{Lefschetz pairs}.
Any pair of this list corresponds to one of the
intersection points, and holds the smallest
and the largest indices of the wires intersected at this point,
numerated locally near the intersection point (see \cite{MoTe1} and \cite{GaTe}).

For example, in the wiring diagram of Figure \ref{wdtolp},

\FIGUREy{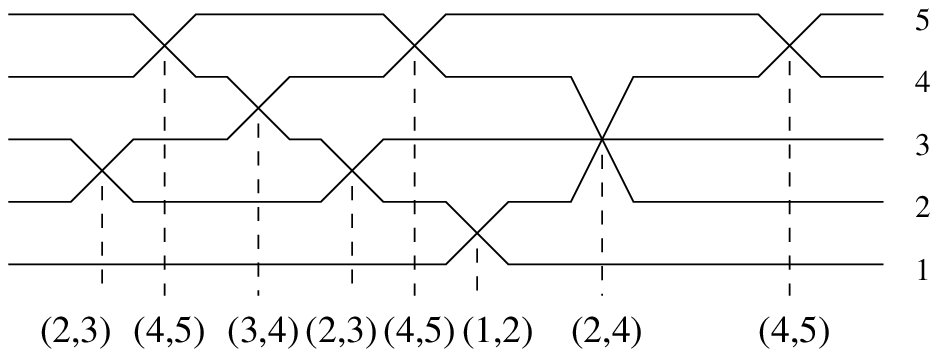}{3cm}{wdtolp}

\noindent
the list of Lefschetz pairs is:
$$( \Lpair{4}{5},\Lpair{2}{4},\Lpair{1}{2},\Lpair{4}{5},\Lpair{2}{3},\Lpair{3}{4},\Lpair{4}{5},\Lpair{2}{3} ).$$
To every list of Lefschetz pairs there is a corresponding wiring
diagram, which is constructed by the reverse procedure.

At some points in this paper it would be necessary to slightly
change the definition of a wiring diagram, by fixing some angles
$\alpha_1 < \dots < \alpha_\ell$, and assuming that the \th{i}
wire is 'broken' at some distant positive point at angle
$\alpha_i$, and at some distant negative point at angle
$-\alpha_i$; all with respect to the guiding line. This is done to
prevent the wires from meeting at infinity; whenever we discuss
the
projective 
 fundamental group, the diagrams are meant to have
this structure.

\subsection{The incidence lattice of an arrangement}
Let $\cL = \set{L_1,\cdots, L_{\ell}}$ be an arrangement of
lines. By ${\rm Lat}(\cL)$ we denote the partially-ordered set of
non-empty intersections of the $L_i$, ordered by inclusion (see
\cite{OT}). We include the whole plane and the empty set in ${\rm
Lat}(\cL)$, so that it becomes a lattice.

For example, if we numerate the intersection points of the
arrangement of lines in Figure \ref{latowd}, as presented in
Figure \ref{latolat},

\FIGUREy{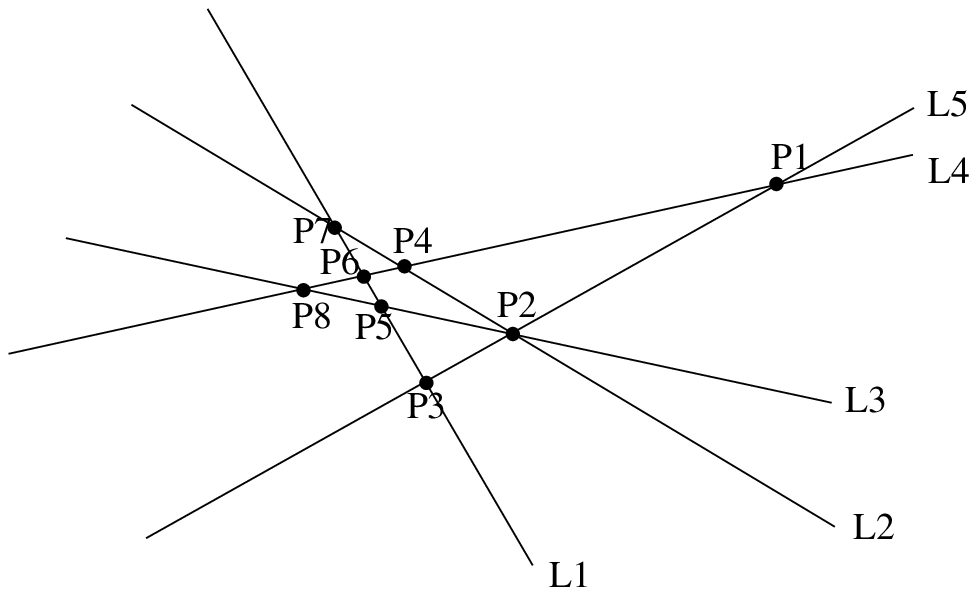}{3cm}{latolat}

\noindent
we get the lattice presented in Figure \ref{figlat}.

\FIGUREy{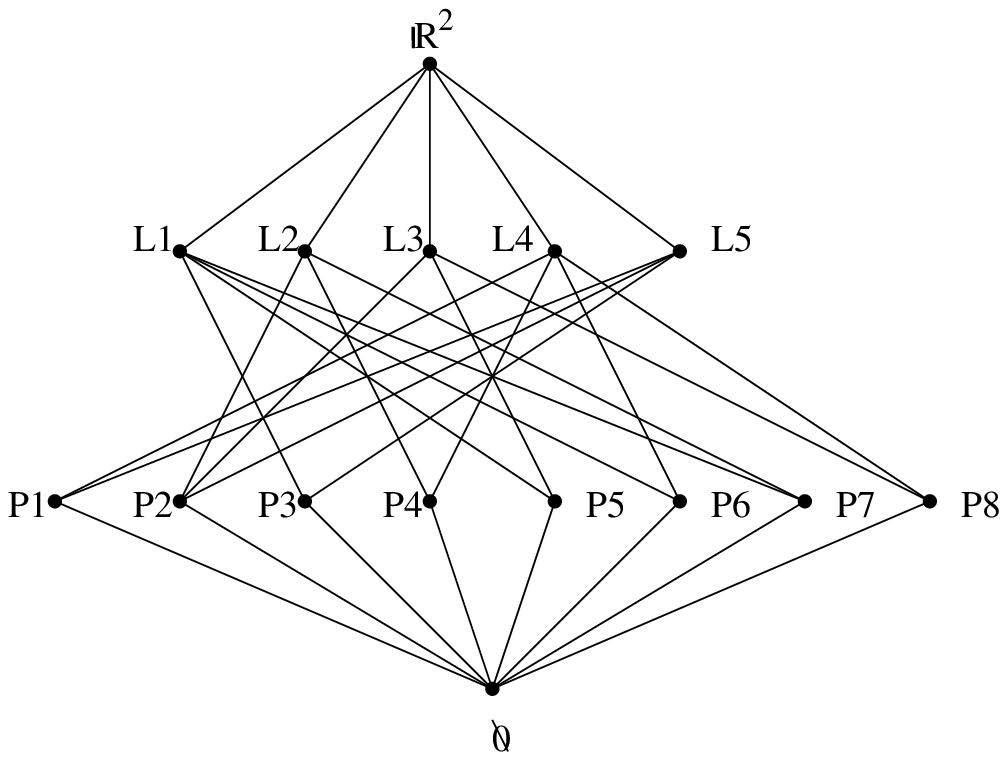}{5cm}{figlat}


\subsection{The signature}\label{sign}

We define another combinatorial invariant of line arrangements: \defin{the
signature}.
\bde
The {\bf signature} of an arrangement is $[2^{n_2}3^{n_3}
\cdots]$ where
$n_k$ is the number of points in which $k$ lines intersect.

We make the agreement to omit every $c^0$.

\ede
Obviously, any wiring diagram induced by an arrangement of lines
has the property that every two wires intersect each other
exactly once. This is the \defin{unique intersection property}.

\bre\label{UIP}
Note that from the unique intersection property it follows
that the order of the lines in $- \infty$ is the inverse of that in $\infty$.
\ere

An easy property for the signature $[2^{n_2}3^{n_3}\dots]$ of a wiring diagram
with $\ell$ lines is that
\begin{equation}\label{SUIP}
\sum_{k\geq 2}n_k{{k} \choose {2}} = {{\ell}\choose {2}}
\end{equation}

\noindent
(there are ${{\ell}\choose {2}}$ simple intersection points in a
generic line arrangement with $\ell$ lines, and every point of
multiplicity $k$ replaces ${{k} \choose {2}}$ simple intersection points).

\section{Computation of the fundamental group}\label{MT}

In this section we present the computation of a the fundamental
groups of the complement of (complexified) wiring diagram. This
is an easy generalization of the Moishezon-Teicher method
\cite{MoTe1} and the van Kampen theorem, which are usually used
to compute the fundamental group of the complement of line
arrangements. The algorithm is used in the proof of Theorems
\ref{si_preserve} and \ref{arel_preserve} below.
A proof for the algorithm itself can be found in
\cite{MoTe1}.

\medskip

Let $D$ be a closed disk in $\R^2$, $K \sbs \Int(D)$ a set of $\ell$ points, and
$u \in\partial D$.
Let $\cB$ be the group of all diffeomorphisms
$\be: D\ra D$ such that $\be|_{\partial D}$ is the identity and $\be(K)=K$.
The action of such $\be$ on the disk applies to paths in $D$, which induces 
an automorphism on 
$\pi_1(D-K,u)$.
The \defin{braid group},
$B_{\ell}[D,K]$, is the group $\cB$ modulo the subgroup of
diffeomorphisms inducing the trivial automorphism on
$\pi_1(D-K,u)$. An element of
$B_{\ell} [D,K]$ is called a \defin{braid}.

Let $\si \sbs \Int(D)$ be a path connecting two points $a,b \in
K$, which avoids all the other points in $K$. Using
$\si$ one defines a diffeomorphism of $D$ by exchanging
$a,b$ along the path (more precisely, along two paths parallel to
$\si$). The resulting braid is called a
\defin{half-twist}, and denoted by $H(\si)$. It can be seen that
$B_{\ell}[D,K]$ is generated by the half-twists. For simplicity,
we will assume that $D = \set{ z \in \C \suchthat \power{z -
\frac{\ell+1}{2}} \leq \frac{\ell+1}{2} }$, and that $K = \set{ 1, 2, \cdots , \ell } \sbs
D$.

Choose a point $u_0 \in D$ (for convenience we choose it to be
below the real line). The group $\pi_1(D-K,u_0)$ is freely
generated by $\Gamma_1,\dots,\Gamma_\ell$, where $\Gamma_i$ is a
loop starting and ending at $u_0$, enveloping the \th{i} point in
$K$. The set $\set{\Gamma_1,\dots,\Gamma_\ell}$ is called a
\defin{geometric base} or \defin{g-base} of $\pi_1(\C^2-K,u_0)$.

Let $(\Lpair{a_1}{b_1},\dots,\Lpair{a_p}{b_p})$ be a list of
\LPs\ associated to a wiring diagram ${\bf w}$ with $\ell$ wires. The
\fg\ of the complement of the diagram is a quotient group of
$\pi_1(D-K,u_0)$. There are $p$ relations, one for every intersection point.
In each point we will compute an object called a
\defin{skeleton}, from which the relation is computed.


Let $\ba = (\Lpair{a_1}{b_1},\cdots,\Lpair{a_p}{b_p})$ be a list of Lefschetz pairs.
In order to compute the skeleton $s_i$ associated to the
\th{i} intersection point $\Lpair{a_i}{b_i}$
we start with an \defin{initial skeleton} which is presented in
Figure \ref{fig3}, in which the points correspond to the lines of
the corresponding wiring diagram.

\FIGUREy{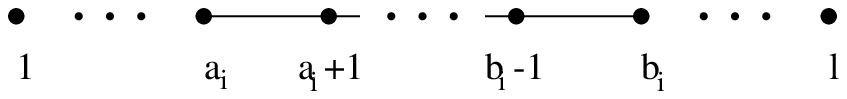}{0.7cm}{fig3}

To this skeleton we apply the Lefschetz pairs
$\Lpair{a_{i-1}}{b_{i-1}},
\cdots,\Lpair{a_1}{b_1}$. A Lefschetz pair  $\Lpair{a_j}{b_j}$
acts by rotating the region from $a_j$ to $b_j$ by $180^{\circ}$
counterclockwise without affecting any other points.

For example, consider the list $\ba =
(\Lpair{2}{3},\Lpair{2}{4},\Lpair{4}{5},\Lpair{1}{3},\Lpair{3}{4})$.
The initial skeleton for $\Lpair{3}{4}$ is given in Figure
\ref{fig4}.

\FIGUREx{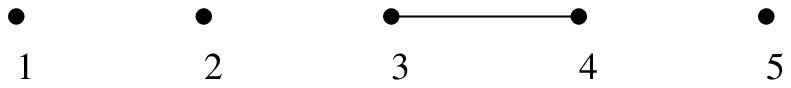}{6cm}{fig4}

Applying $\Lpair{1}{3}$ and then
$\Lpair{4}{5}$, we get the skeleton of Figure \ref{fig5}.

\FIGUREx{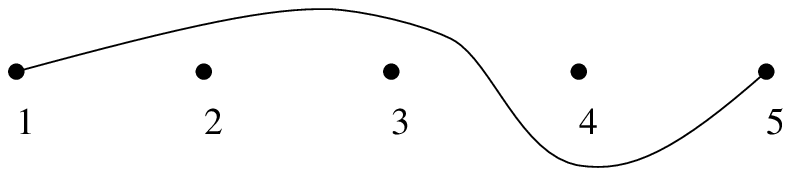}{6cm}{fig5}

Then, applying  $\Lpair{2}{4}$ yields the skeleton of Figure
\ref{fig6}, and finally
acting with  $\Lpair{2}{3}$ we get the skeleton in Figure
\ref{fig8}.

\FIGUREx{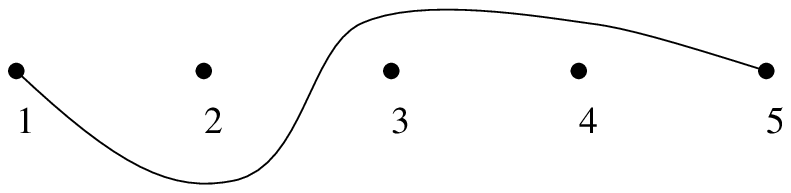}{6cm}{fig6}

\FIGUREx{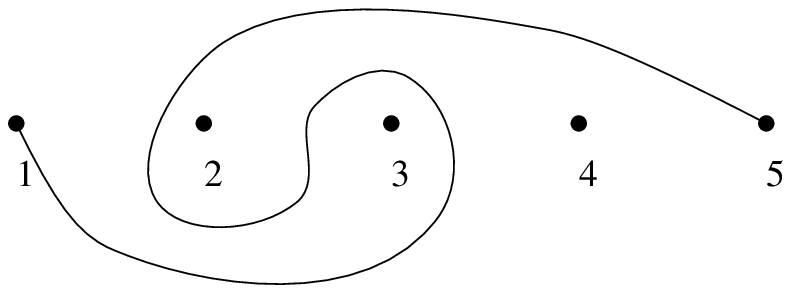}{6cm}{fig8}


From the resulting skeleton we compute the relation, as follows.
We first explain the case when $\Lpair{a_i}{b_i}$ corresponds to a
simple point, \ie\ $b_i - a_i =1$. Then the skeleton is a path
connecting two points.

Choose an arbitrary point on the path and 'pull' it down,
breaking the path into two parts, which are connected in one end
to $u_0$ and in the other to the two end points in $K$.

The loops associated to these two paths are elements in the group
$\pi_1 (D-K,u_0)$, and we call them $a_1$ and $a_2$. These are
elements in $\pi_1(D-K,u_0)$, which commute in the \fg.

Figure \ref{av_bv} illustrates this procedure.

\FIGUREy{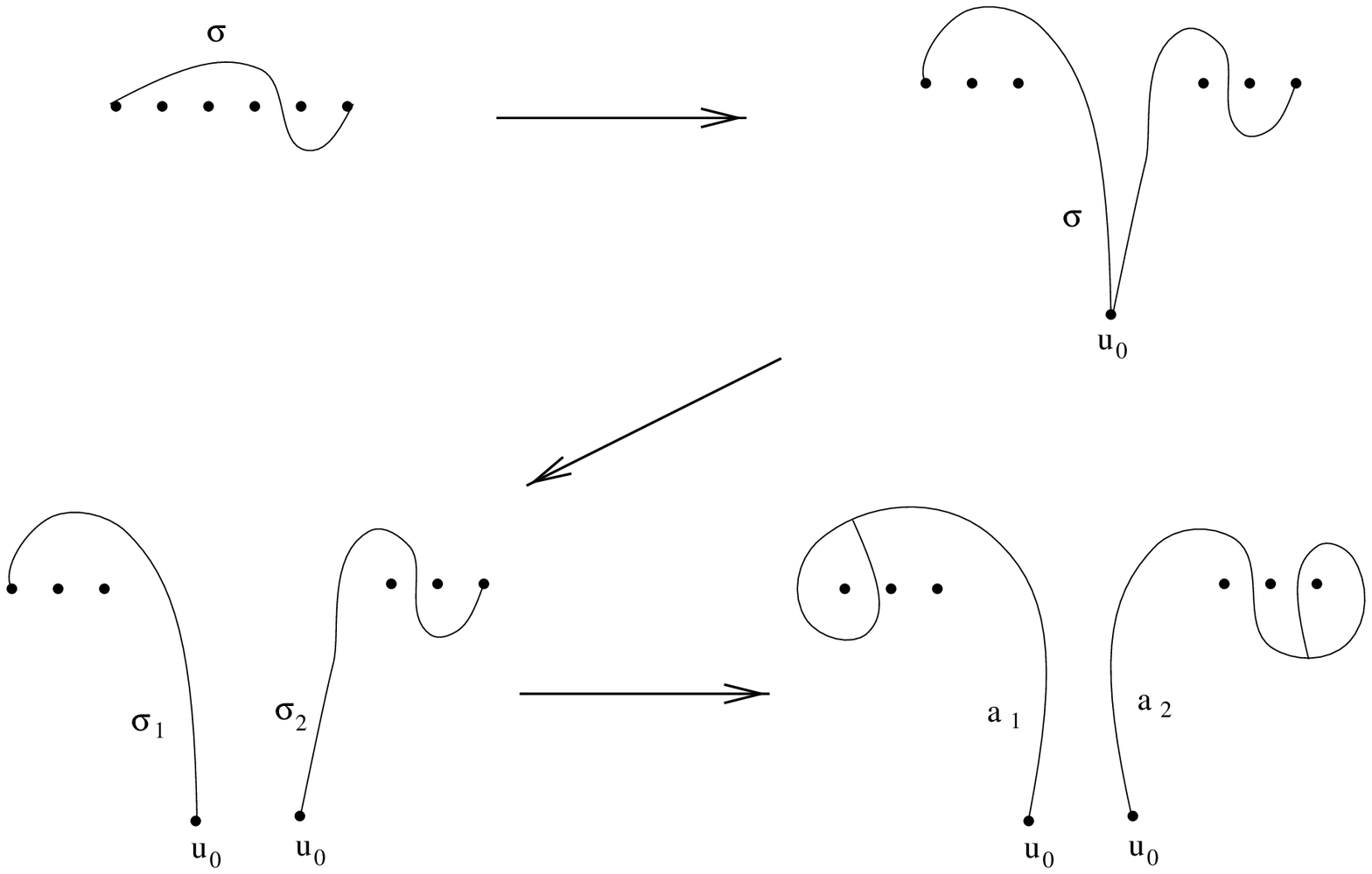}{6cm}{av_bv}

Now we show how to write $a_1$ and $a_2$ as words on the
generators
$\set{\Ga_1, \cdots, \Ga_\ell}$ of $\pi _1(D-K,u_0)$. We begin
with the generator corresponding to the end point of
$a_1$ (or $a_2$), and conjugate it as we move along $a_1$ (or
$a_2$) from its end point on $K$ to $u_0$ as follows: for every
point $i \in K$ which we pass from above, we conjugate by $\Ga_i$
when moving from left to right, and by $\Ga_i^{-1}$ when moving
from right to left.

For example, in the above figure,
$$a_1 = \Ga_3 \Ga_2 \Ga_1 \Ga_2^{-1} \Ga_3^{-1}, \quad a_2 = \Ga_4 ^{-1} \Ga_6 \Ga_4$$
and so the induced relation is:
$$\Ga_3 \Ga_2 \Ga_1 \Ga_2^{-1} \Ga_3^{-1} \cdot \Ga_4 ^{-1} \Ga_6 \Ga_4 =\Ga_4 ^{-1} \Ga_6 \Ga_4 \cdot \Ga_3 \Ga_2 \Ga_1 \Ga_2^{-1} \Ga_3^{-1}$$

One can check that the relation is independent of the point in
which the path is broken.

\bigskip

For a multiple intersection point we compute the elements in the
group
$\pi_1 (D-K,u_0)$ in a similar way, but the induced relations are
 of the type
$$a _k a_{k-1} \cdots  a _1 =
 a _1 a _k \cdots  a _3 a _2 = \cdots =
 a _{k-1} a _{k-2} \cdots a _1 a _k.$$
We choose an
arbitrary point on the path and pull it down to
$u_0$.
For each of the $k$ end points of the skeleton, we generate the
loop associated to the path from
$u_0$ to that point, and translate this path to a word on
$\G_1,\dots,\G_\ell$ by the procedure given above.

\FIGUREy{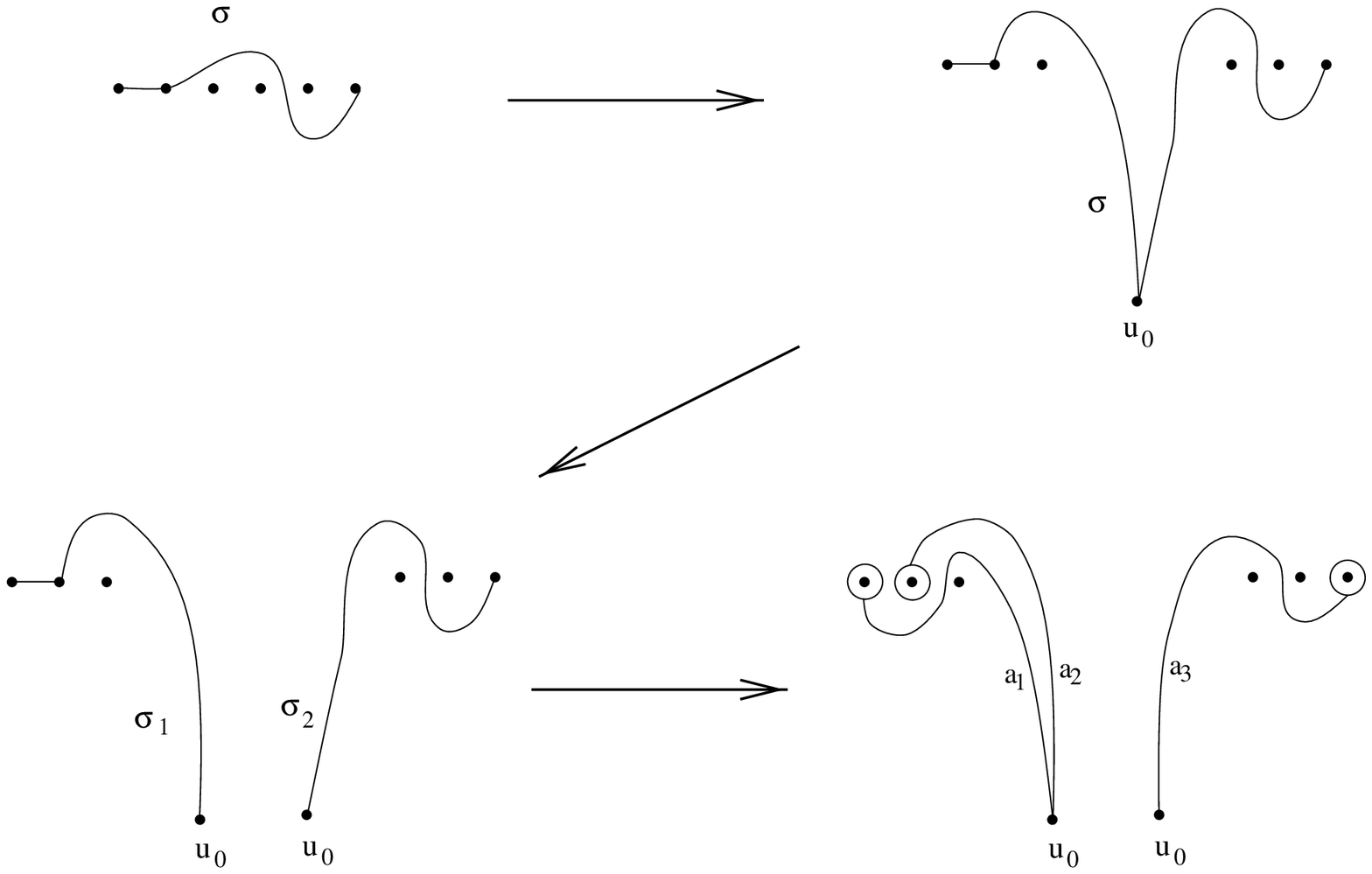}{6cm}{av_bv_mul}

In the example given in Figure \ref{av_bv_mul}, we have
$a_1 =
\G_3 \G_1 \G_3^{-1}$, $a_2 = \G_3 \G_2
\G_3^{-1}$ and $a_3 = \G_4^{-1} \G_6 \G_4$, so the relations are
\begin{eqnarray*}
\G_4^{-1} \G_6 \G_4 \cdot \G_3 \G_2 \G_3^{-1} \cdot \G_3 \G_1
\G_3^{-1} & = & \G_3 \G_1 \G_3^{-1} \cdot \G_4^{-1} \G_6 \G_4 \cdot
\G_3 \G_2 \G_3^{-1}\\
& = & \G_3 \G_2 \G_3^{-1} \cdot \G_3 \G_1
\G_3^{-1} \cdot \G_4^{-1} \G_6 \G_4.
\end{eqnarray*}

Finally, the projective fundamental group $\pi _1 (\C\PP ^2 - {\bf
w})$ is the quotient of the affine fundamental group $\pi _1 (\C
^2 - {\bf w})$ which we just computed, by the  relation
$$\Ga _{\ell} \cdots \Ga _1 =1.$$

\section{Relations and actions on wiring diagrams} \label{actions}

Fix a signature $S=[2^{n_2}3^{n_3}\dots]$. Denote by $W_S$ the set of all lists of
Lefschetz pairs with that given signature, for which the associated
wiring diagram has the unique intersection property (defined in Subsection \ref{sign}).
Note that the number of lines and intersection points is determined by $S$: there are
$p=\sum{n_k}$ points, and the number of lines can be computed from Equation (\ref{SUIP}).

In this section we introduce two equivalence relations and three
actions on $W_S$ or its quotient sets. The motivation for the
equivalence relations and the actions comes from the identities in
the braid group and the isometries of the projective plane.

\subsection{Disjoint intersection points}\label{ss:equiv}

We say that two Lefschetz pairs
$\Lpair{a}{b}$ and $\Lpair{c}{d}$ are \defin{disjoint} if the
correpsonding integral segments $[a,b]$ and $[c,d]$
are disjoint. If two adjacent pairs are disjoint, then the
corresponding intersection points have no common lines.

We say that two lists of Lefschetz pairs ${\bf a, b} \in W_S$ are \defin{equivalent},
and denote ${\bf a} \equiv \bb$,
if it is possible to reach from ${\bf a}$
to $\bb$ by switching adjacent disjoint pairs. This is an equivalence
relation on $W_S$.

For example, the following two wiring diagrams are equivalent:

\FIGUREy{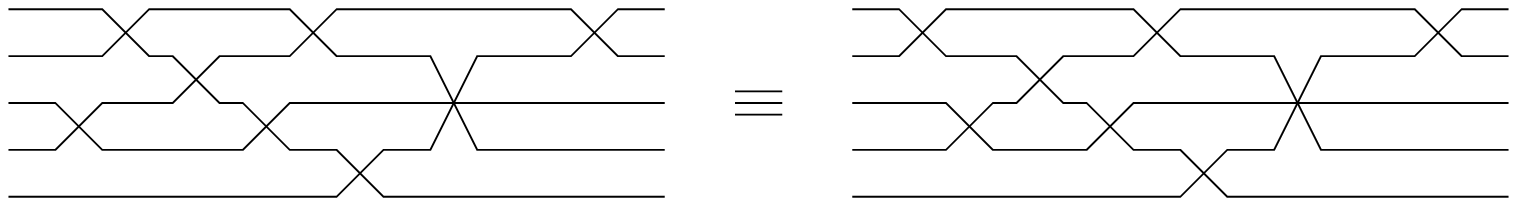}{1.5cm}{figequiv}

%
%

Just as we have $\si_i \si_j = \si_j \si_i$ if $|j-i|>1$ for the standard generators
of Artin's braid group, here too exchanging
two disjoint points does not change the wiring diagram
topologically, so we have
\begin{prop}
If $\ba \equiv \ba'$ then $\ba$ and $\ba'$ have the same incidence
lattice, and the same affine and projective \fgs.
\end{prop}

\subsection{Reflection}\label{ss:tau} 

We define an action $\tau$ on $W_S$, motivated by a reflection of
the plane by a line perpendicular 
to the guiding line:

\bde
For ${\bf a} = (\Lpair{a_1}{b_1},\Lpair{a_2}{b_2}, \cdots , \Lpair{a_p}{b_p}) \in W_S$, set
$$\tau({\bf a}) := (\Lpair{a_p}{b_p},\Lpair{a_{p-1}}{b_{p-1}}, \cdots , \Lpair{a_1}{b_1}).$$
\ede

\bre \label{w_d_tau}
If ${\bf a} \equiv \bb$ then it is easy to see that
$\tau({\bf a}) \equiv \tau(\bb)$, so that $\tau$ is well-defined on
$\Wmodequiv$.
\ere

Obviously $\tau$ does not change the wiring diagram
topologically, so we have
\begin{prop}
The wiring diagrams $\ba$ and $\tau(\ba)$ have the same incidence
lattice, and the same affine and projective \fgs.
\end{prop}

\subsection{Rotation}\label{ss:mu} 

The next action is related to rotation of the plane and its
effect on wiring diagrams. Consider a smooth rotation of a line
arrangement, when the guiding line is kept fixed. When we make a
rotation, the projections of the intersection points slide along
the guiding line. It may happen that two projections coincide
(when two disjoint intersection points are one above the other;
at that moment, the list of Lefschetz pairs induced from the
wiring diagram is not defined). When we continue the rotation,
the two projections switched places and we get a new wiring
diagram, equivalent to the original one (in the sense of
Subsection \ref{ss:equiv}).

A more radical \label{geomu} and interesting change happens when one of the lines
(the one numbered $1$, when the rotation is clockwise) becomes
perpendicular to the guiding line. Immediately before it happens,
the projections of the intersection points on this line are consecutive,
and immediately after the rotation their order is reversed. The local numbers
of the lines from the left to these points are increased by one,
and those from the right are decreased by one.

\medskip

\noindent
{\bf Example.} The arrangement in Figure \ref{lasformu}(a) was rotated clockwise,
until the thick line is perpendicular to the guiding line, and
then some. We get the arrangement given in Figure
\ref{lasformu}(b).

\FIGUREy{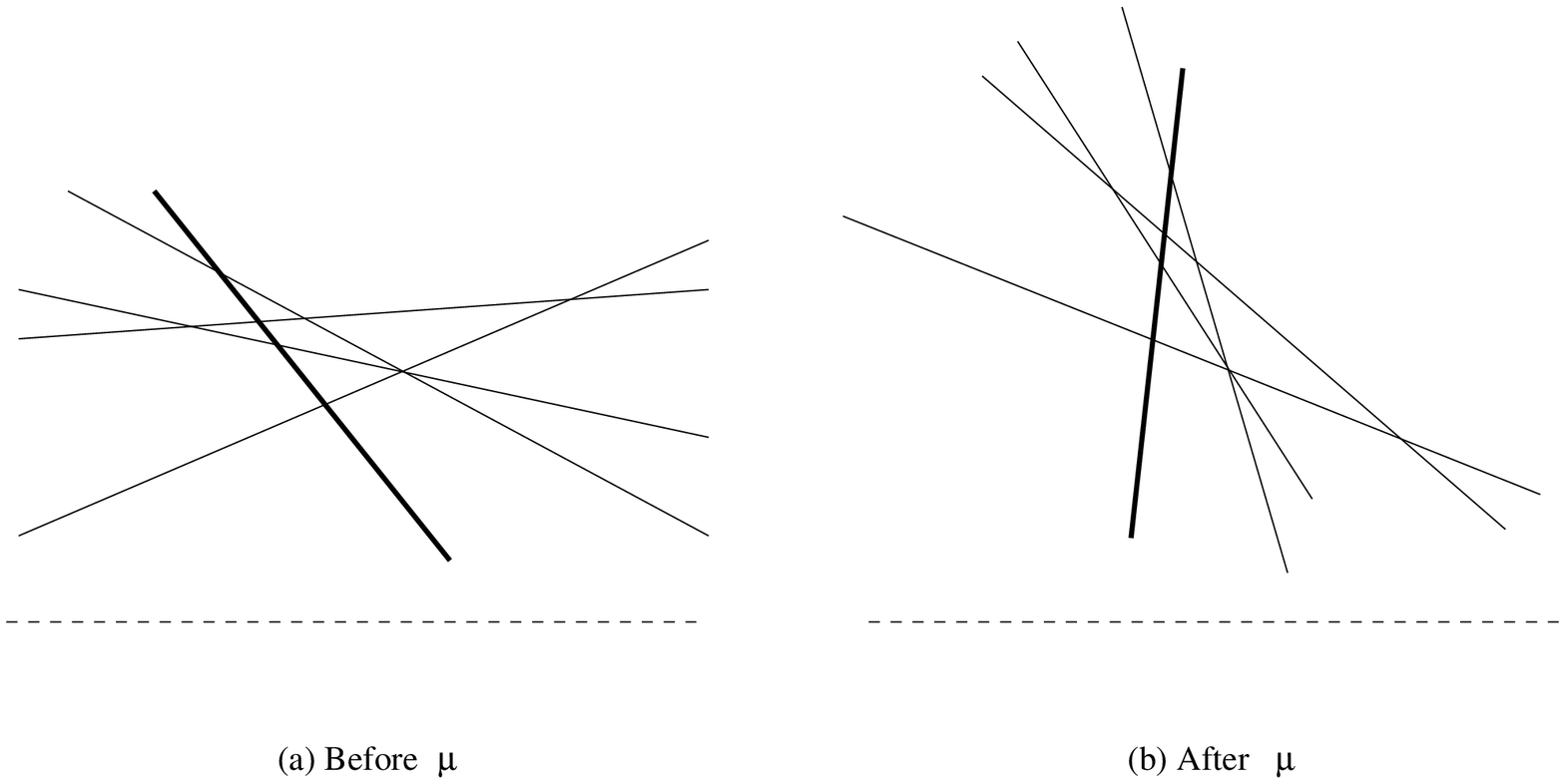}{5cm}{lasformu}

\FIGUREy{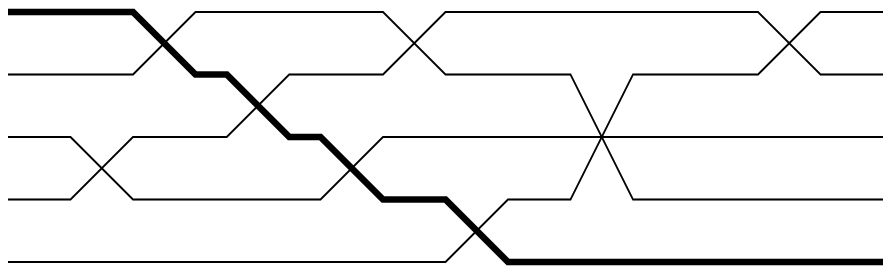}{2.5cm}{wdbeforemu}

The wiring diagram corresponding to Figure \ref{lasformu}(a) is
given in Figure \ref{wdbeforemu}.
 The corresponding list of
Lefschetz pairs is
$$( \Lpair{4}{5},\Lpair{2}{4},\Lpair{1}{2},\Lpair{4}{5},\Lpair{2}{3},\Lpair{3}{4},\Lpair{4}{5},\Lpair{2}{3} ),$$
which can be reordered to the following 
equivalent list:
$$( \Lpair{4}{5},\Lpair{2}{4},\Lpair{4}{5},\Lpair{1}{2},\Lpair{2}{3},\Lpair{3}{4},\Lpair{4}{5},\Lpair{2}{3} ) $$

The wiring diagram corresponding to Figure \ref{lasformu}(b) is
given in Figure \ref{wdaftermu},

\FIGUREy{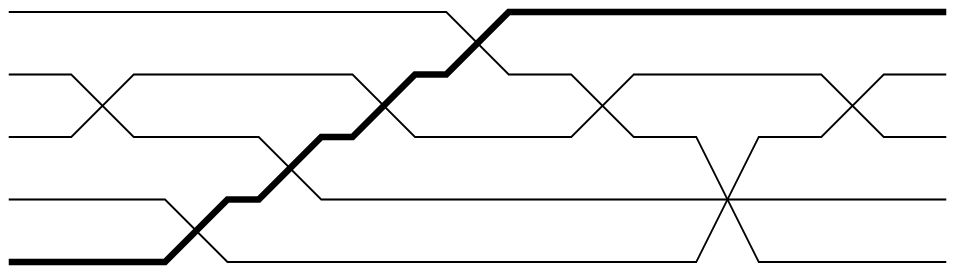}{2.5cm}{wdaftermu}

\noindent
with the corresponding list of Lefschetz pairs being
$$(\Lpair{3}{4},\Lpair{1}{3},\Lpair{3}{4}, \Lpair{4}{5},\Lpair{3}{4},\Lpair{2}{3},\Lpair{1}{2}, \Lpair{3}{4}).$$

\medskip

This  definition for wiring diagrams induced by a line arrangement
will now be generalized to arbitrary wiring diagrams with
the unique intersection property.

Let ${\bf a} = (\Lpair{a_1}{b_1},\Lpair{a_2}{b_2}, \cdots ,
\Lpair{a_p}{b_p})$ be an element of $W_S$. We inductively
construct a decomposition of the list of pairs in ${\bf a}$ into
three disjoint sub-lists $L_+,L_-$ and $L_0$. We view every
$\Lpair{a_i}{b_i}$ as a permutation acting on the indices
$\set{ 1, \cdots, \ell }$, where  $\Lpair{a_i}{b_i}$ sends every $a_i \leq t \leq b_i$ to
$a_i+b_i-t$, and leaves the other indices fixed.

Set $x=1$ and $L_+=L_-=L_0=\emptyset$. For each $i=1, \dots ,p$, act as follows.
If $a_i \leq x \leq b_i$, add $\Lpair{a_i}{b_i}$ to the list $L_0$, and set $x=a_i+b_i-x$.
If $x>b_i$, add $\Lpair{a_i}{b_i}$ to the list $L_-$, and if $x<a_i$, add $\Lpair{a_i}{b_i}$
to the list $L_+$. Continue to the next value of $i$.

During this procedure, the line numbered $1$ at $\infty$ always carries
the local index $x$. As a result, whenever $a_i \leq x \leq b_i$, we actually
have $x=a_i$ (since otherwise the line with local number $a_i$ would have to intersect
the first line at some point before, and again at point number $i$, contradicting the unique intersection property).

If follows that the intersection points which were combined as $L_0$ all lay
on the first line, the points in the list $L_+$ are above that line, and the ones in $L_-$
are below it.

As an example, see Figure \ref{wdbeforemu}, where we have:
$$L_{+} = (\Lpair{4}{5},\Lpair{2}{4},\Lpair{4}{5}) \ ; \ L_0=(\Lpair{1}{2},\Lpair{2}{3},\Lpair{3}{4},\Lpair{4}{5}) \ ; \ L_- = (\Lpair{2}{3})$$

It is now obvious that we can reorder {\bf a} into an equivalent list ${\ba'}$,
which can be written as an \defin{ordered} union ${\bf a}'= L_+ \cup L_0 \cup L_-$.
Thus we can define an action $\mu$ as follows.

\bde\label{def_mu}
Let ${\bf a} \in W_S$ be a list of Lefschetz pairs,
$${\bf a} = (\Lpair{a_1}{b_1},\Lpair{a_2}{b_2}, \cdots , \Lpair{a_p}{b_p}).$$
Decompose ${\bf a} \equiv L_+ \cup L_0 \cup L_-$ as above, so that we can write:
$$L_+ = (\Lpair{a_1}{b_1}, \cdots , \Lpair{a_u}{b_u}),$$
$$L_0 = (\Lpair{a_{u+1}}{b_{u+1}}, \cdots , \Lpair{a_v}{b_v}),$$
$$L_- = (\Lpair{a_{v+1}}{b_{v+1}}, \cdots , \Lpair{a_p}{b_p}).$$
where $1 \leq u < v \leq n$.

We substruct one from the indices in $L_+$, invert $L_0$ and add one to the indices in $L_-$; $\mu({\bf a})$ is defined by:
\begin{eqnarray*}
\mu({\bf a}) & := & (\Lpair{a_1-1}{b_1-1}, \cdots , \Lpair{a_u-1}{b_u-1}, \\
    & & \qquad \Lpair{a_v}{b_v},\Lpair{a_{v-1}}{b_{v-1}}, \cdots ,\Lpair{a_{u+1}}{b_{u+1}}, \\
    & & \qquad \qquad \Lpair{a_{v+1}+1}{b_{v+1}+1}, \cdots , \Lpair{a_p+1}{b_p+1}).
\end{eqnarray*}
\ede

In order for the definition to make sense, we must show that
$\mu ({\bf a})$ is indeed a valid list of Lefschetz pairs (a list
that can be induced from a wiring diagram with the unique
intersection property), with the same signature.

This can be verified as follows. After reordering {\bf a} to its
equivalent form ${\bf a}'= L_+ \cup L_0 \cup L_-$, break the two
rays of line number 1, and glue them back in exactly the opposite
direction. Then rotate the middle section of this line until the
order of the intersection points it carries is reversed. This
transformation results in a new wiring diagram, which is easily
seen to be the one corresponding to $\mu (\ba)$.

\bre \label{w_d_mu}\label{mukeeps}
If ${\bf a}' \equiv {\bf a}$, then $\mu ({\bf a}') = \mu ({\bf a})$.
\ere
\begin{proof}
The three sublists $L_+$, $L_0$ and $L_-$ are the same for every
$\ba' \equiv \ba$. Therefore, $\mu (\ba ')$ is the same as well.
\end{proof}

Since $\mu$ only shifts the numbers of lines, we have
\begin{rem}
The wiring diagrams $\ba$ and $\mu(\ba)$ have the same incidence
lattice.
\end{rem}

Finally, we prove the following:

\bpr
Let $\ba$ be a list of Lefschetz pairs, and let {\bf w} be its
associated wiring diagram. Let ${\bf w_{\mu}}$ be the wiring
diagram associated to the list of Lefschetz pairs $\mu(\ba)$. Then
$$\pi _1 (\C ^2 - {\bf w}) \cong \pi _1 (\C ^2 - {\bf w_{\mu}})$$
and
$$\pi _1 (\C\PP ^2 - {\bf w}) \cong \pi _1 (\C\PP ^2 - {\bf w_{\mu}})$$
\epr

\begin{proof}
Since equivalent wiring diagrams have the same fundamental
groups, we may reorder ${\bf w}$ so that all the intersection
points on the first wire are consecutive (as explain in the
definition of $\mu$).  Then, by definition, the wiring diagram
${\bf w_\mu}$ is generated from $\bw$ by breaking the first wire
of
$\bw$ at the end points of its two rays. We can smoothly rotate
them until we achieve the wiring diagram ${\bf w_\mu}$. Keeping
the same angle of both rays along the rotation, the diagram is
always topologically the same (even in $\C\PP^2$), and the result
follows.
\end{proof}


\subsection{Beyond infinity}\label{ss:si}

Our third action, denoted by $\si$, is well-defined on $W_S$ but
not on $\Wmodequiv$.

Consider a wiring diagram associated to $\ba \in W_S$. Recall that from the unique
intersection property it follows that the order of the pseudolines at $\infty$ is
the inverse of that at $- \infty$.
Glue the wires together by identifying the two
far edges of every wire. This gives an embedding of the wiring diagram on a
M\"obius
strip. To recapture \ba, we simply cut the strip at $\infty$.

The action $\si$ is defined as follows: slide the point $\infty$ to the left
(alternatively, rotate the strip to the right), until it crosses one
intersection point. Then cut the strip at the new place of $\infty$.

The corresponding map on the set of lists of Lefschetz pairs is
as follows (recall that the order of the lines is inverted in
infinity). Let $\ell$ be the number of wires in the wiring
diagram. Denote by $J= \Lpair{1}{\ell}$ the following permutation:
$$J(i)=\ell+1-i, \quad 1 \leq i \leq \ell$$

\bde
Let $\ba = (\Lpair{a_1}{b_1},\Lpair{a_2}{b_2}, \cdots ,
\Lpair{a_p}{b_p}) \in W_S$. Then we define
$$\si(\ba) := (\Lpair{a_2}{b_2}, \cdots ,
\Lpair{a_p}{b_p},\Lpair{J(b_1)}{J(a_1)}).$$
\ede

\medskip

To show that $\si (\ba)$ is again induced by a wiring diagram satisfying the unique
intersection property,
write $\ba = (\pi _1, \cdots, \pi _p)$ where $\pi _i = \Lpair{a_i}{b_i}$
is the permutation that sends every $a_i \leq t \leq b_i$ to
$a_i+b_i-t$, and leaves the other indices fixed.
Note that the signature of $(\pi_1,\dots,\pi_p)$ is given by the
multiset of the sizes of supports of $\pi_i$, and, assuming the
signature satisfies Equation \eq{SUIP}, $(\pi_1,\dots,\pi_p)$
corresponds to a wiring diagram with the unique intersection
property iff
$\pi_1 \dots \pi_\ell = J$. Now,
$\si (\ba) = (\pi _2,
\cdots,
\pi _p, J
\pi_1 J^{-1})$, and

\begin{eqnarray}
\pi _2 \cdots \pi _p \cdot J \pi_1 J^{-1} & = & \pi _1^{-1} (\pi_1 \pi _2 \cdots \pi _p) J \pi_1 J^{-1} =
\nonumber \\
& = & \pi _1^{-1} J^2 \pi _1 J^{-1} = J \nonumber
\end{eqnarray}
since $J^2$ is the identity.

\medskip

\noindent
{\bf Examples.}
Consider the wiring diagram corresponding to
$$\ba = (\Lpair{1}{2},\Lpair{2}{3},\Lpair{1}{2}),$$
with signature $[2^3]$. Then $\si(\ba) =
(\Lpair{2}{3},\Lpair{1}{2},\Lpair{2}{3})$ and
$\si ^2
(\ba)=\ba$.

\medskip

Here is a slightly more complicated example. In the following pair
of line arrangements, the right one is obtained from the other by
pushing the rightmost point ``through infinity'' until the two
bold lines are parallel, then letting them intersect at the left
side of the line arrangement.

\FIGUREy{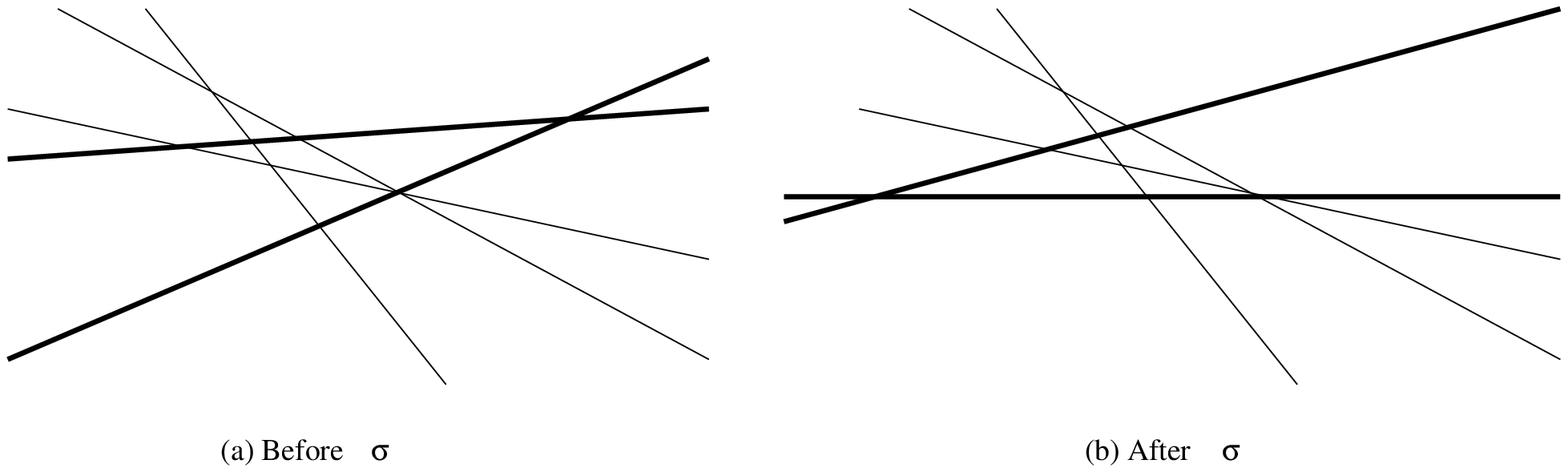}{3.6cm}{lasforsi}

\medskip

The \wid\ corresponding to \ref{lasforsi}(a) is given in Figure
\ref{wdbeforesi}.
\FIGUREy{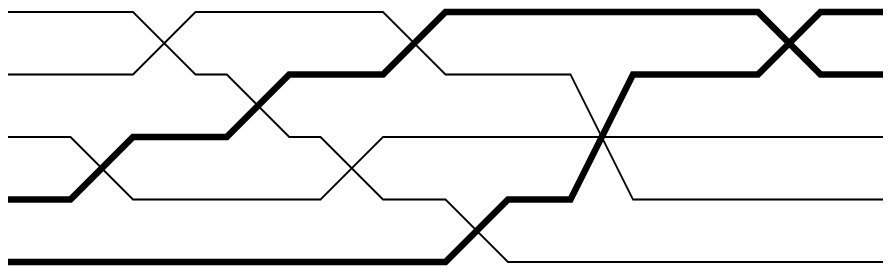}{2.4cm}{wdbeforesi}
If we apply the above described  action on the wiring diagram, we
get the diagram in Figure \ref{wdaftersi}.
\FIGUREy{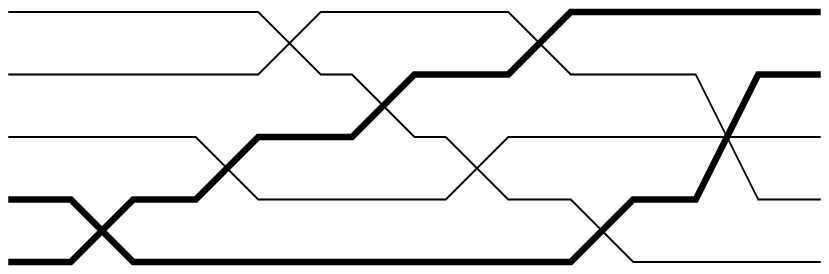}{2.5cm}{wdaftersi}
These figures demonstrate a simple point going 'through
$\infty$', but the same procedure applies to a multiple
intersection point as well.

As in previous actions, acting with $\si$ changes only the
numeration of the points (and may also change the numbers of the
lines, if the rightmost point involves the lowest line), so we
have
\begin{rem}
For every $\ba \in W_S$, $\ba$ and $\si(\ba)$ have the same
incidence lattice.
\end{rem}

Next, we prove that $\si$ preserves the fundamental groups.

\bth \label{si_preserve}
For every wiring diagram $\ba \in W_S$, $\ba$ and $\si(\ba)$ have
the same affine and projective fundamental groups.
\eth

\begin{proof}
We will in fact give an explicit isomorphism between the
fundamental groups which were described in Section \ref{MT}.

Let
$\ba = ( \Lpair{a_1}{b_1}, \cdots,\Lpair{a_p}{b_p} )$ be a list of Lefschetz
pairs, and let ${\bf w}$ and ${\bf w_{\si}}$ be the wiring
diagrams associated to $\ba$ and
$\si (\ba)$, respectively. Fix $a,s$ such that $\Lpair{a}{a+s} =
\Lpair{a_1}{b_1}$.

\forget
When we apply $\si$, its geometric effect is exchanging the
elements of this g-base which correspond to the wires intersected
at the rightmost intersection point (this change behaves very
similar to a generalization of Hurewitz move on this g-base, see
\cite{MoTe1}). This geometric effect induces an algebraic map
between the two fundamental groups $\pi_1(\C ^2 -{\bf w})$ and
$\pi_1(\C ^2 -{\bf w_{\si}})$ as follows.

According to van Kampen theorem,
we associate one generator to each  wire in the diagram. Then,
for all the wires not intersected at the point corresponds to the
Lefschetz pair $\Lpair{a_1}{b_1}$, their associated generators in
$\pi _1(\C^2 -{\bf w})$ will be mapped to the corresponding
generators in $\pi _1(\C^2 -{\bf w_{\si}})$. For the wires which
are intersected at that point, their associated generators in
$\pi _1(\C^2 -{\bf w})$ will be mapped to the elements in
$\pi _1(\C^2 -{\bf w_{\si}})$
which obtained by an Hurewitz move on the g-base of $\pi _1(\C^2
-{\bf w})$.

Without loss of generality, we may assume that $\Lpair{a_1}{b_1}
= \Lpair{a}{a+s}$, then the appropriate generalized Hurewitz move
is:


Therefore, the induced map $\phi: \pi _1 (\C ^2-{\bf w}) \to \pi
_1 (\C ^2-{\bf w_{\si}})$ is:
\forgotten

Denote by $\Gamma_1,\dots,\Gamma_\ell$ the geometric generators of
$\pi_1(\C^2-\bw)$, and by $\Gamma_1',\dots,\Gamma_\ell'$ the
generators of $\pi_1(\C^2-\bw_\sigma)$. We define
$\phi:\pi_1(\C^2-\bw) \ra \pi_1(\C^2-\bw_\sigma)$ by
\begin{center}
$$\phi (\Ga _i) =  \left\{ \begin{array}{cl}
{\Ga _i' }          &  { 1 \leq i < a } \\
{{\Ga' _a} ^{-1} \cdots {\Ga'_{a+s-1}}^{-1} \Ga' _{a+s} \Ga'_{a+s-1} \cdots  \Ga_a'}  & {i=a} \\
{{\Ga' _a} ^{-1} \cdots {\Ga'_{a+s-2}}^{-1} \Ga_{a+s-1}' \Ga_{a+s-2}' \cdots \Ga _a'} & {i=a+1} \\
{\cdots} & {\cdots} \\
{{\Ga' _a} ^{-1} \Ga_{a+1}' \Ga _a'} & {i=a+s-1} \\
{\Ga _a' } & {i=a+s} \\
{\Ga _i' }          &  {  a+s < i \leq \ell }
\end{array}
\right.
$$
\end{center}

This definition is motivated by Figure \ref{hurewitz}, which
corresponds to an appropriate Hurwitz move \cite{Hu}.

\FIGUREy{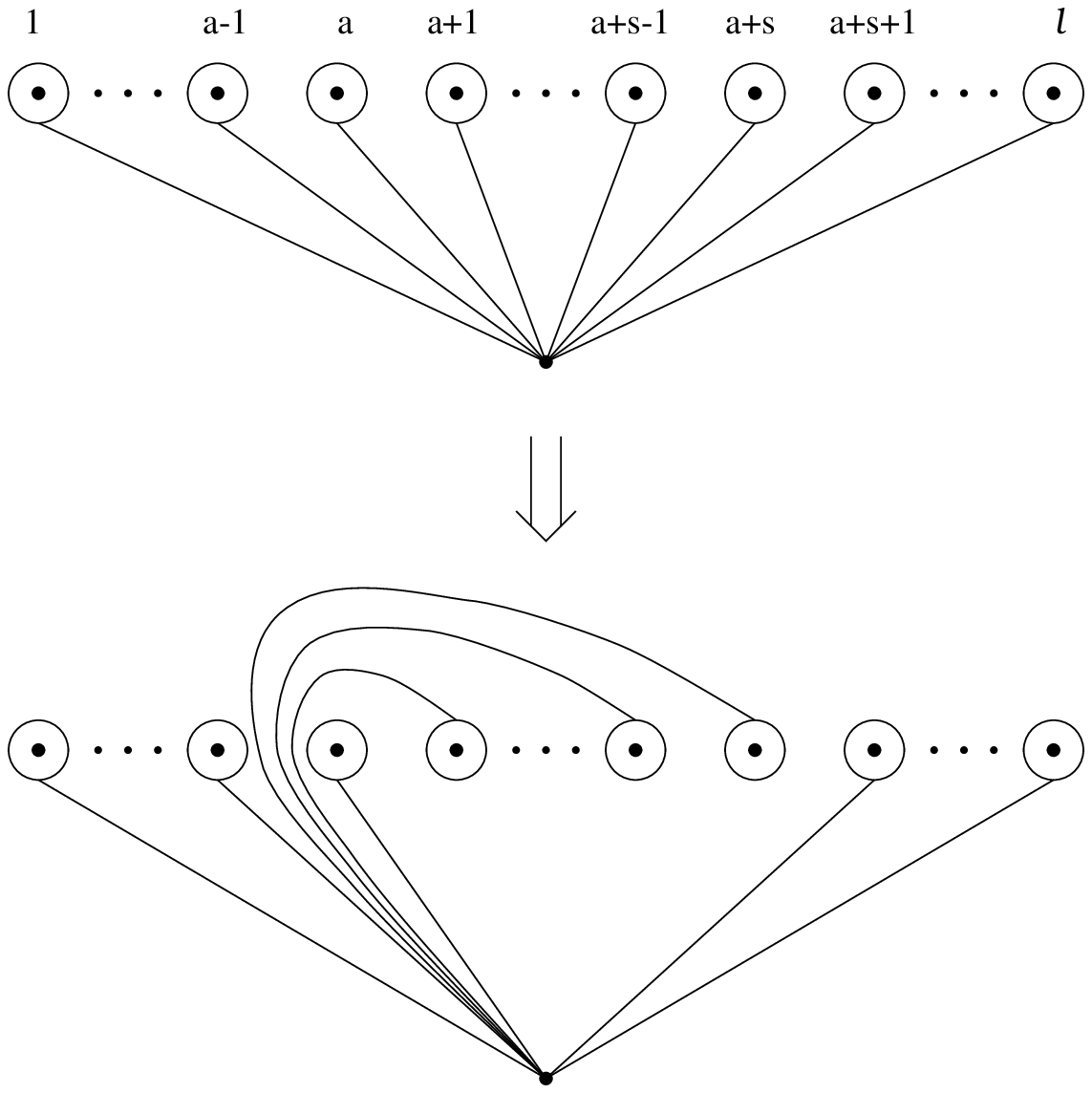}{6cm}{hurewitz}

In order to show that $\phi$ is well defined and is an
isomorphism, we need to show that it carries the set of relations
of $\pi_1(\C^2-\bw)$ to the relations of
$\pi_1(\C^2-\bw_\sigma)$.

Recall that by definition
$$\si(\ba) = (\Lpair{a_2}{b_2}, \cdots , \Lpair{a_p}{b_p},\Lpair{J(b_1)}{J(a_1)}).$$

\forget
Using the algorithm in section \ref{MT}, we compute a finite
presentation for the groups $\pi _1(\C ^2-{\bf w})$ and $\pi
_1(\C ^2-{\bf w_{\si}})$ and we show that the above canonical map
is indeed an isomorphism map between the two groups.

\medskip

Let $\set{ \Gamma_1,\Gamma_2,\dots,\Gamma_\ell }$ and
$\set{ \Gamma_1',\Gamma_2',\dots,\Gamma_\ell' }$ be the g-bases of
$\pi_1(D-K,*)$ and $\pi_1(D'-K',*)$ respectively, which their images are,
by the van Kampen theorem, the sets of generators of the (affine
and projective) \fgs\ of
${\bf w}$ and ${\bf w_{\si}}$ respectively.

Now we have to show that this canonical map is an isomorphism
between the two groups. Because the map $\phi$ is a bijection
between the two sets of generators, it is sufficient to show that
the relations in
$\pi _1(\C ^2- {\bf w})$  are carried to the relations in $\pi _1(\C ^2- {\bf w_{\si}})$.
\forgotten

We will first show that the relations in $\pi _1(\C^2 -{\bf w})$
associated to the points
$\Lpair{a_2}{b_2},\dots,\Lpair{a_p}{b_p}$ are equivalent to the relations
in $\pi _1(\C^2 -{\bf w_{\si}})$ associated to the these points.

Let $2 \leq j \leq p$. According to the algorithm in section
\ref{MT}, the \th{j} skeleton of $\ba$ is obtained by applying on the initial skeleton
$\Lpair{a_j}{b_j}$ the
halftwists corresponding to $\Lpair{a_{j-1}}{b_{j-1}}$, then to
$\Lpair{a_{j-2}}{b_{j-2}}$, and so on, down to $\Lpair{a_1}{b_1}$.
At the same time, the \th{(j-1)} skeleton of
${\bf w_\si}$ is obtained from the same initial skeleton,
by applying $\Lpair{a_{j-1}}{b_{j-1}}$ down to
$\Lpair{a_2}{b_2}$, without applying $\Lpair{a_1}{b_1}$.
Therefore, in order to get the same skeleton in both cases, one
has to act on the \th{j} skeleton of
{\bf w} by the inverse to the halftwist corresponds to
$\Lpair{a_1}{b_1}$. The action of the Hurwitz move (and the 
map $\phi$) corresponds
to this inverse action: it turns the local region in the disk
clockwise (in the direction opposite to that of the halftwist's
action). Therefore one gets the same skeletons, and hence
equivalent relations in the presentation of the fundamental
groups.

It remains to show that the relation induced by the skeleton
corresponding to the first point, $\Lpair{a_1}{b_1}$, 
in {\bf w}, is equivalent to the relation induced by the skeleton corresponds 
to the last point, $\Lpair{J(a_1)}{J(b_1)}$, 
of ${\bf w_\si}$.

According to the algorithm, the skeleton corresponding to
$\Lpair{a_1}{b_1}$ for {\bf w} is simply the straight segment between $a_1$ and
$b_1$. The induced
relation is (see
\cite{GaTe}):
$$\Ga _{b_1} \cdots \Ga _{a_1} = \Ga _{b_1-1} \cdots \Ga _{a_1} \Ga _{b_1} = \cdots = \Ga _{a_1} \Ga _{b_1} \cdots \Ga _{a_1+1}.$$

For ${\bf w_{\si}}$, the skeleton is obtained by applying
halftwists corresponding to the rest of the pairs, on the initial
skeleton
$\Lpair{J(b_1)}{J(a_1)}$.

In the braid group, the action of all the $p$ halftwists of a list
of Lefschetz pairs (with the unique intersection property) is
equal to the action of the general halftwist $H([1,\ell])$, so
the skeleton in this case can be obtained as follows: begin with
the initial skeleton
$\Lpair{J(b_1)}{J(a_1)}$, apply $H(\Lpair{1}{\ell})$ and then apply
the inverse of the halftwist corresponding to
$\Lpair{a_1}{b_1}$. Applying $H(\Lpair{1}{\ell})$ on
$\Lpair{J(b_1)}{J(a_1)}$ yields the skeleton corresponding to the
segment $\Lpair{a_1}{b_1}$. Then, applying the halftwist
corresponds to $\Lpair{a_1}{b_1}$ will give the same result.
Therefore we get the same skeleton, and hence equivalent sets of
relations. This proves that the two affine groups are isomorphic.

Finally, compute that
$\phi(\Gamma_\ell\dots\Gamma_1) = \Gamma_\ell'\dots\Gamma_1'$, so
the relation defining the projective group of $\bw$ transfers to
the corresponding relation for $\bw_\sigma$.
\end{proof}


\subsection{Going above and below}

We define another relation on $W_S$, which can be motivated by
the triple relation $\si _i \si _{i+1} \si _i = \si _{i+1} \si _i
\si _{i+1}$ in the braid group.
This relation corresponds to a line passing over a simple point. 
Surprisingly, it turns out that this change of the geometry of
the arrangement, does not change the incidence lattice and
fundamental groups.
%
%

We generalize this to a line passing over a multiple point. Let
$t\geq 2$ denote the number of lines intersecting in a point,
where the lines are locally numerated $c,\dots,c+t-1$ (this point
will be called here the \defin{central} point). Now suppose
another line crosses these
$t$ lines, with all the intersection points adjacent to the
central point. Let
$c+i$ denote the index of this new line before the intersections,
as in Figure
\ref{mul_diag}. We claim that it
does not matter if the new line crosses above or below that
central point.

\FIGUREy{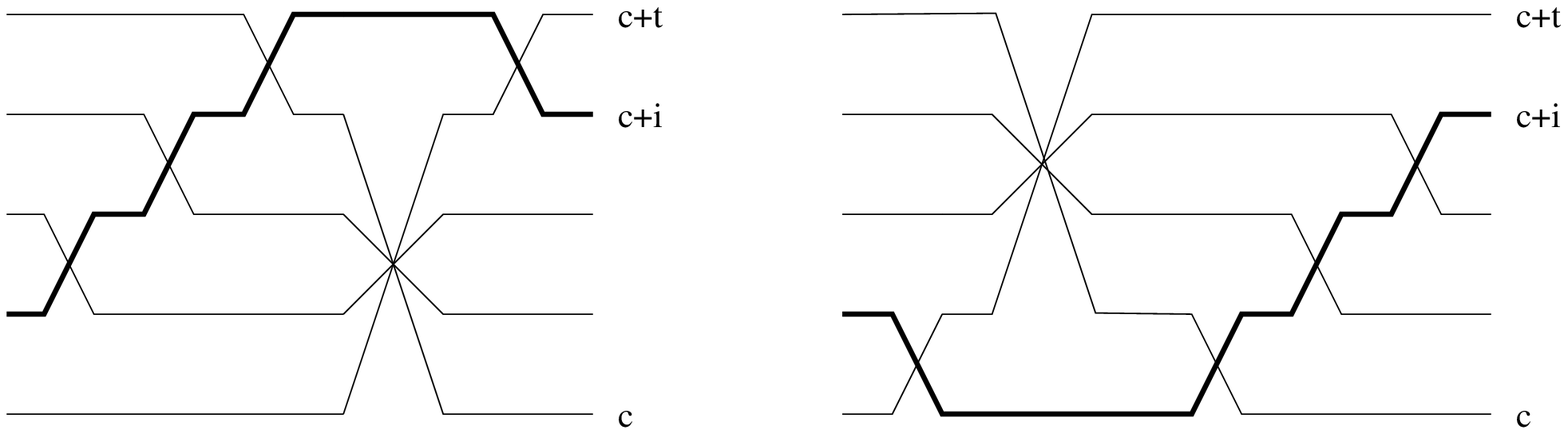}{2.8cm}{mul_diag}

In order to simplify the notation, we write
$c+\Lpair{a}{b}$ instead of $\Lpair{c+a}{c+b}$, and more
generally if $\ba$ is a list of \LPs, by $c+\ba$ we mean adding
$c$ to all the pairs in
$\ba$.
In the situation described above we have
$1 \leq c \leq c+i \leq c+t \leq \ell$, and the lists of
\LPs\ are $\tru{i}{t}$ when the new line crosses from above, and
$\trd{i}{t}$ when it crosses from below, where we define
\begin{eqnarray*}
\tru[]{i}{t} & = & (\Lpair{i}{i+1}, \cdots, \Lpair{t-1}{t},\Lpair{0}{t-1},\\
& & \Lpair{t-1}{t}, \cdots, \Lpair{t-i}{t-i+1}),\\
\trd[]{i}{t} & = & (\Lpair{i-1}{i},\cdots,\Lpair{0}{1},
\Lpair{1}{t},\\
& & \Lpair{0}{1}, \cdots,\Lpair{t-1-i}{t-i}).
\end{eqnarray*}

\bde
Let ${\bf a, b}\in W_S$. We say that ${\bf a} \arel {\bf b}$, if
there is a chain of replacements of
$\tru{i}{t}$ by $\trd{i}{t}$, or \viceversa, which goes from $\ba$ to $\bb$.
\ede


Since this relation changes only the positions of some
intersection points, it obviously preserves the incidence lattice.
We claim that $\arel$ also preserves the affine and projective
fundamental groups.

\begin{thm}\label{arel_preserve}
Let ${\ba,\bb} \in W_S$ be two lists of Lefschetz pairs  such that
$\ba \arel \bb$. Let $\bw_{\ba}$ and $\bw_{\bb}$ are
the wiring diagrams associated to ${\bf a}$ and ${\bf b}$
respectively. Then  ${\bf w_a}$ and ${\bf w_b}$ have the same
affine and projective fundamental groups.
\end{thm}

\begin{proof}
Let ${\bf a}$ and
${\bf b}$
be two lists of Lefschetz pairs, such that
${\bf a} \arel {\bf b}$.
We have to show that
$\pi _1 (\C ^2 -{\bf w_a})\cong \pi _1 (\C ^2 -{\bf w_b})$.
It is obviously enough to assume the only difference between $\ba$
and
$\bb$ is one replacement of $\tru{t}{i}$ by $\trd{t}{i}$, and by
applying $\sigma$ we may assume that this replacement happens in
the first $t+1$ pairs of each list. Then we have
$$\ba = (\tru{i}{t}, \Lpair{a_{t+2}}{b_{t+2}},\cdots,\Lpair{a_p}{b_p}),$$
and
$$\bb = (\trd{i}{t}, \Lpair{a_{t+2}}{b_{t+2}},\cdots,\Lpair{a_p}{b_p}).$$

We will compute the presentations of $\pi _1 (\C ^2 -{\bf w_a})$
and $\pi _1 (\C ^2 -{\bf w_b})$ using the algorithm given in
Section \ref{MT}, where $\set{
\Gamma_1,\Gamma_2,\dots,\Gamma_\ell }$ are the generators for the
first group, and $\set{ \Gamma_1',\Gamma_2',\dots,\Gamma_\ell'}$
for the second. We will show that the map $\Gamma_j \mapsto
\Gamma'_j$ is an isomorphism between the two groups.

We first show that the relations associated to the intersection
points $\Lpair{a_j}{b_j}$ for $t+1<j\leq \ell$ are the same for
both diagrams. According to the algorithm in section
\ref{MT}, the \th{j} skeleton of ${\bf w_a}$ is obtained by
applying on the initial skeleton, associated to the pair
$\Lpair{a_j}{b_j}$, the composition of the following halftwists:
first, those associated to the pairs
$\Lpair{a_k}{b_k}, k=j-1, \dots, t+2$, and then the $t+1$ halftwists associated to the pairs
in $\tru{i}{t}$. At the same time, the \th{j} skeleton of
${\bf w_b}$ is obtained by applying on the same initial skeleton the same $j-(t+1)$
halftwists, and then the $t+1$ halftwists associated to the pairs
in
$\trd{i}{t}$.

Since in the braid group the composition of the halftwists
associated to the pairs in $\tru{i}{t}$ is equal to the
composition of the halftwists associated to the pairs in
$\trd{i}{t}$ (they both equal to the general halftwist
$H([c,c+t])$), we get the same final \th{j}
skeleton for ${\bf w_a}$ and for ${\bf w_b}$. In particular, we
get the same induced relations in both cases.

\medskip

It remains to show that the relations induced by the first $t+1$
intersection points of ${\bf w_a}$ are equivalent to the
relations induced by the first $t+1$ intersection points of ${\bf
w_b}$. We will compute the relations in both cases.

\begin{prop}
The relations induced on the generators
$\Gamma_1,\dots,\Gamma_\ell$ by the points listed in $\tru{i}{t}$
are
\begin{equation}\label{relw1}
[\Ga _{c+i},\Ga_{c+k}] = 1, \qquad 0 \leq k \leq t
\end{equation}
\begin{eqnarray}\label{relw1c}
& & \Ga_{c+t} \Cdots \Ga_{c} = \\
& & \quad = \Ga_{c+t-1} \Cdots \Ga_{c} \cdot \Ga_{c+t}= \nonumber \\
& & \quad = \cdots = \Ga_{c} \cdot \Ga_{c+t} \Cdots \Ga_{c+1}.
\nonumber
\end{eqnarray}
\end{prop}
\begin{proof}
There are three types of skeletons to consider, corresponding to
point before the central point, points after the central point,
and the central point itself.

1. First let
$1
\leq j
\leq t-i$, so that the point is $c+\Lpair{i+j-1}{i+j}$.
We need to apply the halftwists corresponding to the pairs
$c+\Lpair{i+k-1}{i+k}$, $k=j-1,\dots,1$, on the
initial skeleton which corresponds to
$c+\Lpair{i+j-1}{i+j}$. The resulting skeleton is given in
Figure \ref{w1_fig_case1}.

\FIGUREx{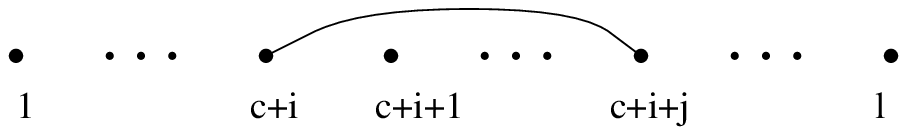}{8cm}{w1_fig_case1}

This skeleton induces the following relation:
$$\Ga _{c+i} \cdot \Ga _{c+i+1}^{-1} \cdots \Ga _{c+i+j-1}^{-1} \Ga _{c+i+j}\Ga _{c+i+j-1} \cdots \Ga _{c+i+1} = $$
$$=\Ga _{c+i+1}^{-1} \cdots \Ga _{c+i+j-1}^{-1} \Ga _{c+i+j}\Ga _{c+i+j-1} \cdots \Ga _{c+i+1} \cdot \Ga _{c+i}.$$

Induction on $j$ now proves Equation
\eq{relw1} for $i < k
\leq t$.

2. We compute the relations for the other points, and then get
back to the central point. Let $t-i+2 \leq j \leq t+1$, and
consider the point
$c+\Lpair{2t-i-j+1}{2t-i-j+2}$ which comes after the central point.
To obtain the skeleton, we begin with the initial skeleton, given
in Figure \ref{w1_fig_case3-1}. We apply the halftwists
corresponding to the pairs
$c+\Lpair{2t-i-k+1}{2t-i-k+2}$, $k=j-1, \cdots, t-i+2$, obtaining the skeleton given in
Figure \ref{w1_fig_case3-2}; then the halftwist corresponding to
the central point
$\Lpair{c}{c+t-1}$, obtaining the skeleton in Figure \ref{w1_fig_case3-3},
and finally the halftwists corresponding to the pairs
$c+\Lpair{i+(k-1)}{i+k}$, $k=t-i, \cdots, 1$, obtaining the final skeleton of Figure
\ref{w1_fig_case3-4}.

\FIGUREx{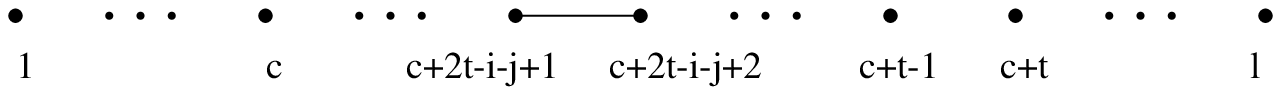}{10cm}{w1_fig_case3-1}

\FIGUREx{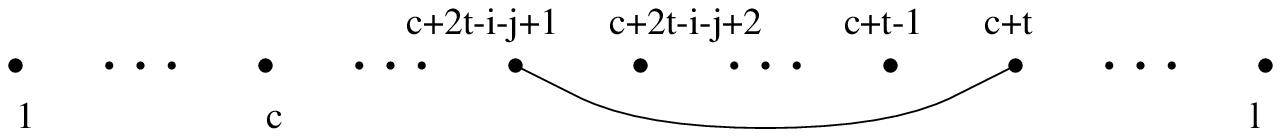}{10cm}{w1_fig_case3-2}

\FIGUREx{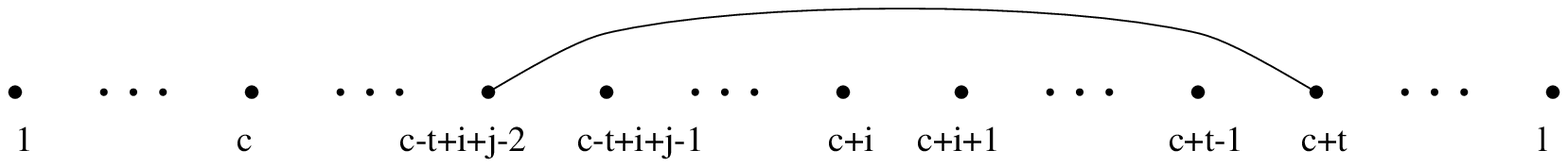}{10cm}{w1_fig_case3-3}

\FIGUREx{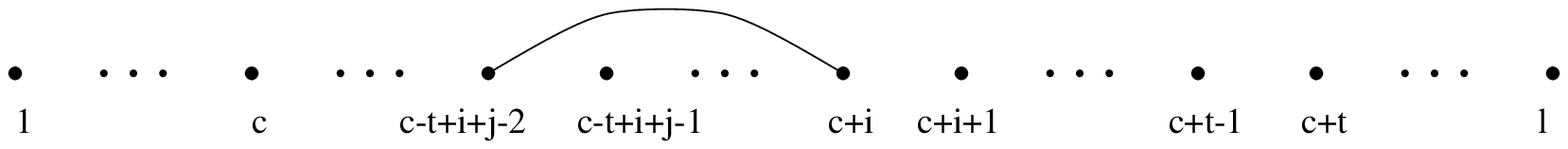}{10cm}{w1_fig_case3-4}

The relation induced by this skeleton is
$$\Ga _{c+i} \cdot \Ga _{c+i-1} \cdots \Ga_{c-t+i+j-1} \Ga_{c-t+i+j-2} \Ga_{c-t+i+j-1}^{-1} \cdots
\Ga _{c+i-1} ^{-1} = $$
$$\Ga _{c+i-1} \cdots \Ga_{c-t+i+j-1} \Ga_{c-t+i+j-2} \Ga_{c-t+i+j-1}^{-1} \cdots
\Ga _{c+i-1} ^{-1} \cdot \Ga _{c+i}.$$

Again by induction on $j$, we obtain the relations
\eq{relw1} for $0 \leq k < i$. Since $\eq{relw1}$ trivially holds
for $k = i$, this set of relations is proved.

3. The central skeleton is obtained by applying the halftwists
correspond to the pairs
$c+\Lpair{i+k-1}{i+k}$,
$k=t-i,\dots,1$ on the initial skeleton corresponding to
$\Lpair{c}{c+t-1}$, which is given in Figure
\ref{w1_fig_case2-1}. The resulting skeleton is
given in Figure \ref{w1_fig_case2-2}.

\FIGUREx{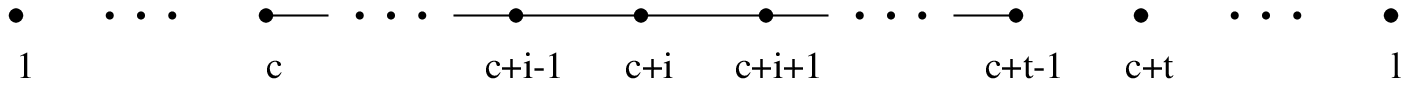}{10cm}{w1_fig_case2-1}

\FIGUREx{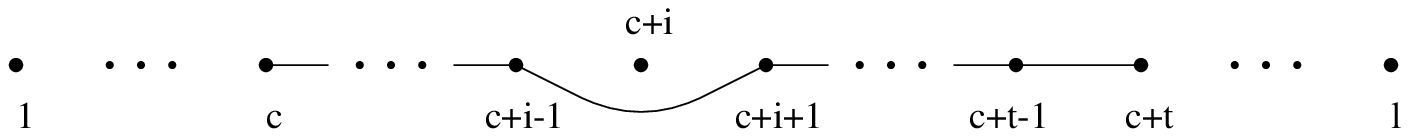}{10cm}{w1_fig_case2-2}

The relation induced by this skeleton is
$$\Ga_{c+t} \cdot \dots \cdot \Ga_{c+i+1}\cdot \Ga_{c+i-1} \cdot \dots \cdot \Ga_{c} = $$
$$= \Ga_{c+t-1} \cdot \dots \cdot \Ga_{c+i+1} \cdot \Ga_{c+i-1} \cdot \dots \cdot \Ga_{c} \cdot \Ga_{c+t}=$$
$$ = \cdots = \Ga_{c} \cdot \Ga_{c+t} \cdot \dots \cdot \Ga_{c+i+1} \cdot \Ga_{c+i-1} \cdot \dots \cdot \Ga_{c+1}.$$
which is Equation \eq{relw1c} since by \eq{relw1} $\Gamma_{c+i}$
commutes with all the generators appearing in this relation.
\end{proof}

Next, we compute the corresponding relations induced by the
skeletons correspond to the first $t+1$ intersection points of
${\bf w_b}$.

\begin{prop}
The relations induced on the generators
$\Gamma_1,\dots,\Gamma_\ell$ by the points listed in $\trd{i}{t}$
are
\begin{equation}\label{relw2}
[\Ga _{c+i}',\Ga_{c+k}'] = 1, \qquad 0 \leq k \leq t
\end{equation}
\begin{eqnarray}\label{relw2c}
& & \Ga_{c+t}' \Cdots \Ga_{c}' = \\
& & \quad = \Ga_{c+t-1}' \Cdots \Ga_{c}' \cdot \Ga_{c+t}' = \nonumber \\
& & \quad = \cdots = \Ga_{c}' \cdot \Ga_{c+t}' \Cdots \Ga_{c+1}'.
\nonumber
\end{eqnarray}
\end{prop}
\begin{proof}
Again we treat the three classes of points separately.

1. Let $1 \leq j \leq i$. The \th{j} skeleton of ${\bf w_b}$ is
obtained by applying the halftwists corresponding to the pairs
$c+\Lpair{i-k}{i-k+1}$, $k=j-1,\dots,1$ on the initial skeleton
induced by the pair $c+\Lpair{i-j}{i-j+1}$. The resulting
skeleton is given in Figure \ref{w2_fig_case1}.

\FIGUREx{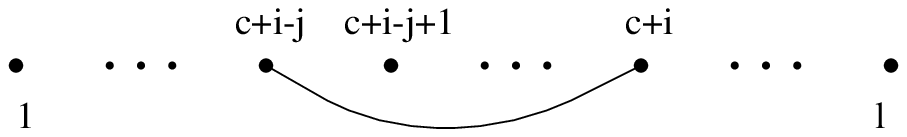}{7cm}{w2_fig_case1}

This skeleton induces the relation
$\Ga' _{c+i-j} \Ga' _{c+i} = \Ga' _{c+i} \Ga' _{c+i-j}$,
and these are the cases $0\leq k < i$ of Equation
\eq{relw2}.

2. Let $i+2 \leq j \leq t+1$. For the \th{j} skeleton we begin
with the initial skeleton corresponding to
$c+\Lpair{i-j}{i-j+1}$ (see Figure
\ref{w2_fig_case3-1}), apply
the halftwists correspond to the pairs
$c+\Lpair{i-k}{i-k+1}$, $k=j-1, \cdots, i-1$ to get the skeleton
of Figure \ref{w2_fig_case3-2}, then apply the central halftwist,
corresponding to $\Lpair{c+1}{c+t}$, and get the skeleton of
Figure \ref{w2_fig_case3-3}, and finally apply the halftwists
correspond to the pairs
$c+\Lpair{i-k}{i-k+1}$, $k=i-1, \cdots, 1$, and get the final
skeleton, presented in Figure \ref{w2_fig_case3-4}.

\FIGUREx{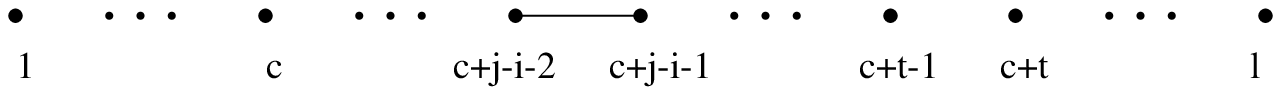}{10cm}{w2_fig_case3-1}

\FIGUREx{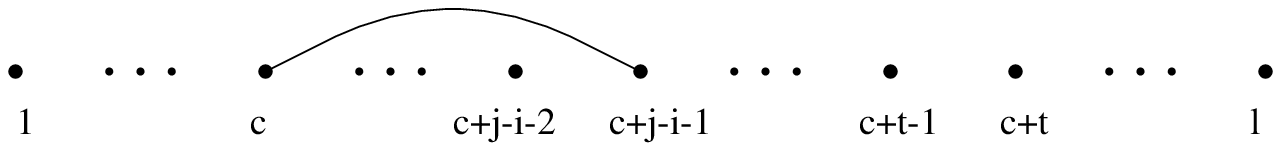}{10cm}{w2_fig_case3-2}

\FIGUREx{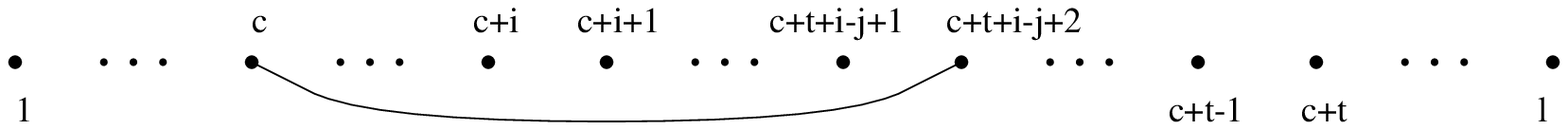}{10cm}{w2_fig_case3-3}

\FIGUREx{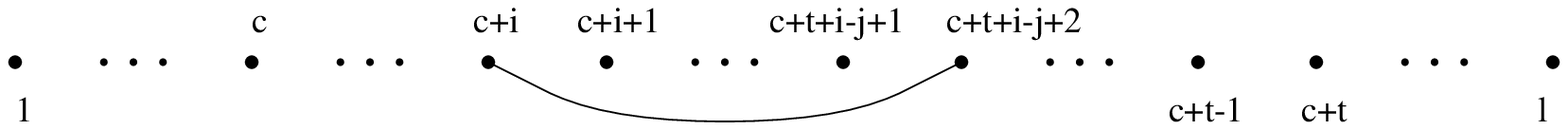}{10cm}{w2_fig_case3-4}

This skeleton induces the relation
$$\Ga' _{c+i} \Ga' _{c+t+i-j+2} = \Ga' _{c+t+i-j+2} \Ga' _{c+i},$$
and ranging over the possible values of $j$ we obtain the
relation of Equation \eq{relw2} for $i < k \leq t$.

3. The central skeleton is obtained by applying the halftwists
correspond to the pairs $c+\Lpair{i-k}{i-k+1}$,
$k=i-1,\dots,1$ to the initial skeleton induced by the pair
$\Lpair{c+1}{c+t}$, which is given in Figure
\ref{w2_fig_case2-1}.

\FIGUREx{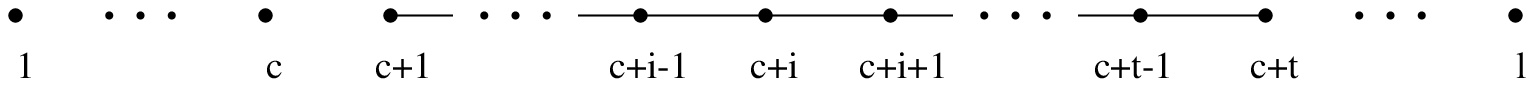}{10cm}{w2_fig_case2-1}

The resulting skeleton is given in Figure
\ref{w2_fig_case2-2}.

\FIGUREx{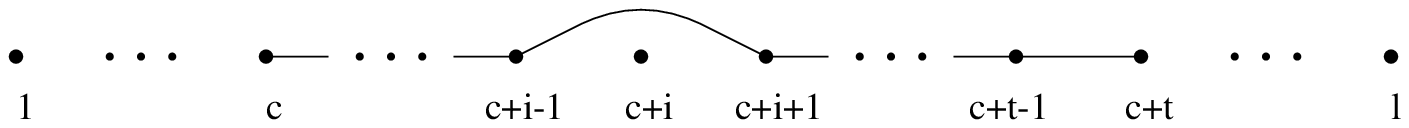}{10cm}{w2_fig_case2-2}  

The corresponding induced relations are
\begin{eqnarray*}
& & \\
& & {\Ga'_{c+i}}^{-1} \Ga'_{c+t} \Ga'_{c+i} \cdot \dots
\cdot {\Ga'_{c+i}}^{-1} \Ga'_{c+i+1} \Ga'_{c+i} \cdot \Ga'_{c+i-1}
\cdot \dots \cdot \Ga'_{c} = \nonumber \\
& & \quad = {\Ga'_{c+i}}^{-1} \Ga'_{c+t-1} \Ga'_{c+i} \cdot \dots
\cdot {\Ga'_{c+i}}^{-1} \Ga'_{c+i+1} \Ga'_{c+i} \cdot \Ga'_{c+i-1}
\cdot \dots \cdot \Ga'_{c} \cdot {\Ga'_{c+i}}^{-1} \Ga'_{c+t}
\Ga'_{c+i}= \nonumber \\
& & \quad = \cdots = \Ga'_{c} \cdot {\Ga'_{c+i}}^{-1}
\Ga'_{c+t}
\Ga'_{c+i} \cdot \dots \cdot {\Ga'_{c+i}}^{-1} \Ga'_{c+i+1}
\Ga'_{c+i} \cdot \Ga'_{c+i-1} \cdot \dots \cdot
\Ga'_{c+1}.\nonumber
\end{eqnarray*}

Using the relations \eq{relw2} we can simplify this to obtain
relation \eq{relw2c}.
\end{proof} 

Comparing the two propositions we discover the same set of relations, so
the affine groups are shown to be isomorphic.

Since our isomorphism is $\Gamma_j\mapsto \Gamma_j'$, the
projective relation
$\Ga_{\ell} \cdots \Ga_1 = 1$ goes to $\Ga'_{\ell} \cdots \Ga'_1 =
1$, so we have also proved that
$\pi _1(\C\PP^2 -{\bf w_a}) \cong \pi _1(\C\PP^2 -{\bf w_b})$.
\end{proof}

\section{Connections between the actions}\label{connections}

In this section we point out two group actions induced by
$\sigma,\tau$ and $\mu$. Note that it is meaningless to ask what
is the group generated by $\sigma$ and $\mu$, since
$\mu$ is defined only up to equivalence, and $\sigma$ is not well defined on equivalence classes.
Let
$D_n$ denote the dihedral group of order $2n$.

\ble\label{D-2l}
Let $S$ be a signature on $\ell$ pseudolines. Then $\mu$ and
$\tau$ induce an action of the dihedral group $D_{2 \ell}$ on the set
$\Wmodequiv$. In particular $\mu ^{2
\ell} =1$,
$\tau ^2 =1$, and
$\tau \mu \tau = \mu ^{-1}$ on this set.
\ele

\begin{proof}
$\tau$ inverts the order of the pairs, so that $\tau ^2$ is the identity even on $W_S$.

By the geometric definition, each application of $\mu$ consists of inverting the direction
of one pseudoline. $\mu ^{\ell} (\ba)$ is a wiring diagram which is
equivalent to $\ba$ rotated by $180 ^{\circ}$.
In particular, $\mu ^{2 \ell}$ is the identity on $\Wmodequiv$.

It remains to check that $(\tau \mu)^2 = {\rm Id}$. Let $\ba$ be
a list of Lefschetz pairs. As in the definition of $\mu$, we can
decompose some $\ba'\equiv \ba$ into three disjoint consecutive
sublists: $\ba'= L_+ \cup L_0 \cup L_-$, where we write:
$L_+ = (\Lpair{a_1}{b_1}, \cdots, \Lpair{a_m}{b_m})$,
$L_0 = (\Lpair{1}{c_1},\Lpair{c_1}{c_2}, \cdots,\Lpair{c_{k-2}}{c_{k-1}}, \Lpair{c_{k-1}}{\ell})$, and
$L_- = (\Lpair{d_1}{e_1}, \cdots, \Lpair{d_p}{e_p})$.

Now, compute:
\begin{eqnarray*}
\tau \mu \tau \mu(\ba') & = & \tau \mu \tau \mu(\Lpair{a_1}{b_1}, \cdots, \Lpair{a_m}{b_m},  \\
    & & \quad \Lpair{1}{c_1},\Lpair{c_1}{c_2}, \cdots,\Lpair{c_{k-2}}{c_{k-1}},\Lpair{c_{k-1}}{\ell}, \\
    & & \quad \Lpair{d_1}{e_1}, \cdots, \Lpair{d_p}{e_p}) =  \\
& = & \tau \mu \tau(\Lpair{a_1-1}{b_1-1}, \cdots, \Lpair{a_m-1}{b_m-1},  \\
    & & \quad \Lpair{c_{k-1}}{\ell},\Lpair{c_{k-2}}{c_{k-1}}, \cdots ,\Lpair{c_1}{c_2},\Lpair{1}{c_1},  \\
    & & \quad  \Lpair{d_1+1}{e_1+1}, \cdots, \Lpair{d_p+1}{e_p+1}) =  \\
& = & \tau \mu(\Lpair{d_p+1}{e_p+1}, \cdots, \Lpair{d_1+1}{e_1+1}, \\
    & & \quad \Lpair{1}{c_1},\Lpair{c_1}{c_2}, \cdots,\Lpair{c_{k-2}}{c_{k-1}}, \Lpair{c_{k-1}}{\ell},  \\
    & & \quad \Lpair{a_m-1}{b_m-1}, \cdots, \Lpair{a_1-1}{b_1-1}) =  \\
 & = & \tau(\Lpair{d_p}{e_p}, \cdots, \Lpair{d_1}{e_1},  \\
    & & \quad \Lpair{c_{k-1}}{\ell},\Lpair{c_{k-2}}{c_{k-1}}, \cdots ,\Lpair{c_1}{c_2},\Lpair{1}{c_1}, \\
    & & \quad \Lpair{a_m}{b_m}, \cdots, \Lpair{a_1}{b_1}) = \\
& = & (\Lpair{a_1}{b_1}, \cdots, \Lpair{a_m}{b_m}, \\
    & & \quad \Lpair{1}{c_1},\Lpair{c_1}{c_2}, \cdots,\Lpair{c_{k-2}}{c_{k-1}}, \Lpair{c_{k-1}}{\ell}, \\
    & & \quad \Lpair{d_1}{e_1}, \cdots, \Lpair{d_p}{e_p}) = \ba'.
\end{eqnarray*}

\end{proof}

\medskip

\ble\label{D-2p}
Assume that the signature $S$ has $p$ intersection points. Then
$\si$ and $\tau$ induce an action of the dihedral group $D_{2p}$ on $W_S$.
In particular $\si^{2p}=1$, $\tau ^2=1$ and $\tau \si \tau ^{-1}=\si^{-1}$.
\ele

\medskip

\begin{proof}
Let $\ba = (\Lpair{a_1}{b_1},\Lpair{a_2}{b_2}, \cdots , \Lpair{a_p}{b_p})$.
By definition,
$$\si^p(\ba) =(\Lpair{J(b_1)}{J(a_1)},\Lpair{J(b_2)}{J(a_2)}, \cdots , \Lpair{J(b_p)}{J(a_p)})$$
so that $\si ^{2p}=1$.

As was shown in the last lemma, the order of $\tau$ is 2.
Therefore, in order to show that $\sg{\si,\tau}$ is dihedral 
of order $4p$, it remains to verify that
$\si \tau \si \tau = {\rm Id}$. Indeed,
\begin{eqnarray}
\si \tau \si \tau (\ba) & = & \si \tau \si \tau(\Lpair{a_1}{b_1},\Lpair{a_2}{b_2}, \cdots , \Lpair{a_p}{b_p}) \nonumber \\
& = & \si \tau \si(\Lpair{a_p}{b_p},\Lpair{a_{p-1}}{b_{p-1}}, \cdots, \Lpair{a_1}{b_1}) \nonumber \\
& = & \si \tau(\Lpair{a_{p-1}}{b_{p-1}}, \Lpair{a_{p-2}}{b_{p-2}}, \cdots, \Lpair{a_1}{b_1}, \Lpair{J(b_p)}{J(a_p)}) \nonumber \\
& = & \si(\Lpair{J(b_p)}{J(a_p)},\Lpair{a_1}{b_1}, \cdots , \Lpair{a_{p-2}}{b_{p-2}},\Lpair{a_{p-1}}{b_{p-1}}) \nonumber \\
& = & (\Lpair{a_1}{b_1}, \cdots , \Lpair{a_{p-2}}{b_{p-2}},\Lpair{a_{p-1}}{b_{p-1}},\Lpair{JJ(a_p)}{JJ(b_p)}) \nonumber \\
& = & (\Lpair{a_1}{b_1}, \cdots , \Lpair{a_{p-2}}{b_{p-2}},\Lpair{a_{p-1}}{b_{p-1}}, \Lpair{a_p}{b_p}) = \ba .\nonumber
\end{eqnarray}

\end{proof}

\medskip

Actually, $\sg{\si}$ can have orbits of size smaller than $2p$,
as shown by the following example.

\medskip

{\bf Example.}
Let $S$ be the signature $[2^3 3^6]$ of an arrangement of
 $7$ pseudolines and $9$ intersection points (Kelly-Moser's configuration
\cite{KM}).
Let
$$\ba = (\Lpair{3}{5},\Lpair{1}{3},\Lpair{5}{6},\Lpair{3}{5},\Lpair{5}{7},\Lpair{2}{3},\Lpair{3}{5},\Lpair{1}{3},\Lpair{5}{6}).$$
One can check that $\si ^6 (\ba)=\ba$. This diagram has orbit of
size $12$ under the action of $\sg{\sigma,\tau}$.

\medskip

In the proof of Lemma \ref{D-2p} we have seen that
$\si ^p$ is a reflection with respect to a line parallel to the
rays of our pseudolines. In particular $\sigma^p$ is well defined
on $\Wmodequiv$ (even though
$\sigma$ is not).
Similarly, we have shown in Lemma \ref{D-2l} that
$\mu ^{\ell}$ corresponds
to a rotation of the diagram by $180^{\circ}$, so their
composition is a vertical reflection which corresponds to $\tau$.
\begin{cor}
As actions of the set $\Wmodequiv$, we have that
$$\tau= \si^p \mu^{\ell}.$$
\end{cor}

As a summary, we now list the number of classes under various sets of actions
and relations. Consider the signature $S = [2^{13} 3^3 4^1]$ (on $8$ lines).

The table below gives the number of equivalence classes of
$\Wmodequiv$ under the various relations. For example, the forth line tells us
that if two wiring diagrams are considered similar whenever there a chain
of moving from a digram to an equivalent diagram (under $\equiv$ or $\arel$),
or of application of $\mu$, then there are $6104$ similarity classes.

\begin{center}
\begin{tabular}{|cccc||c|}
            \hline
             $\sigma$
                & $\tau$
                & $\mu$
                & $\arel$
                & Number of classes \\
            \hline
- & - & - & - & 354880   \\
- & - & - & + & 114379   \\
- & - & + & - & 22180  \\
- & - & + & + & 6104    \\
- & + & - & - & 177440    \\
- & + & - & + & 54539    \\
- & + & + & - & 11090  \\
- & + & + & + & 3076   \\
+ & - & - & - & 5060 \\
+ & - & - & + & 772  \\
+ & - & + & - & 116  \\
+ & - & + & + & 22    \\
+ & + & - & - & 2558  \\
+ & + & - & + & 398  \\
+ & + & + & - & 116  \\
+ & + & + & + & 22  \\
\hline
\end{tabular}
\end{center}

We also computed the incidence lattices and the fundamental
groups for each of the $22$ similarity classes (when all the
relations and actions are considered). There are $5$ different
lattices, and the classes which fall under the same lattice turn
out to have the same fundamental groups.

\begin{\bib}{10}

\bibitem[CS]{CS} Cohen,~ D.~C. and Suciu,~A.~I., {\it The braid monodromy of
   plane algebraic curves and hyperplane arrangements}, Comment. Math.
   Helvetici {\bf 72}, 285-315 (1997).
\bibitem[GaTe]{GaTe} Garber,~D. and Teicher,~M., {\it The fundamental group's
  structure of the complement of some configurations of real line
  arrangements}, Complex Analysis and Algebraic Geometry,
  edited by T.~Peternell and F.-O.~Schreyer, de Gruyter, 173-223 (2000).
\bibitem[GTV]{GTV} Garber,~D., Teicher,~M. and Vishne,~U.,
   {\it $\pi _1$-classification of real arrangements with up to eight lines}, submitted.
\bibitem[Go]{Go} Goodman,~ J.~E., {\it Proof of a conjecture of Burr, Gr\"unbaum and
   Sloane}, Discrete Math. {\bf 32}, 27-35 (1980).
\bibitem[GP]{GP} Goodman,~J.~E. and Pollack,~R., {\it Allowable sequences
   and ordered types in discrete and computational geometry}, in: New trends in
   discrete and computational geometry, edited by J.~Pach, Springer-Verlag,
   103-134 (1993).
\bibitem[Gr]{Gr} Gr\"unbaum,~B., {\it Arrangements and spreads}, Amer. Math.
   Soc., Providence (1972).
\bibitem[Hu]{Hu} Hurwitz,~A., {\it \"Uber Riemann'sche Fl\"achen mit gegebenen verzweigungspunkten},
   Math. Ann. {\bf 39}, 1-60 (1891).
\bibitem[KM]{KM} Kelly,~L.~M. and Moser,~W.~O.~J., {\it On the number of ordinary lines
determined by $n$ points}, Canad. J. Math. {\bf 10}, 210-219 (1958).
\bibitem[Le]{Le} Levi,~F., {\it Die Teillung der projektiven ebene durch gerade oder
  pseudogerade}, Der. Math.-Phys. Kl. S\"achs Akad. Wiss. {\bf 78}, 256-267
  (1926).
\bibitem[MoTe]{MoTe1} Moishezon,~B.~G. and  Teicher,~M., {\it Braid Group
   Technique in Complex Geometry I, Line Arrangements in $\C\PP ^2$},
   Contemporary Math. {\bf 78}, 425-555 (1988).
\bibitem[OrT]{OT} Orlik,~P. and Terao,~H., {\it Arrangements of Hyperplanes},
   Grundlehren {\bf 300}, Springer-Verlag (1992).
\bibitem[Ri]{Ri} Ringel,~G., {\it Teilungen der ebene durch geraden oder topologische
  geraden}, Math. Zeitschrift {\bf 64}, 79-102 (1956).
\end{\bib}

\end{document}